\documentclass[a4paper,11pt]{amsart}

\usepackage[T1]{fontenc}
\usepackage[utf8]{inputenc}
\usepackage{amsmath,amssymb,amsthm,mathtools}
\usepackage{enumitem}
\usepackage{cite}
\usepackage{microtype}
\usepackage{xcolor}
\usepackage[top=2.7cm,bottom=2.7cm,left=2.5cm,right=2.5cm]{geometry}
\usepackage{tikz}
\usetikzlibrary{arrows.meta,positioning}
\usetikzlibrary{arrows.meta,automata,positioning}

\DeclareMathOperator{\Aut}{Aut}
\DeclareMathOperator{\Stab}{Stab}
\DeclareMathOperator{\Core}{Core}
\DeclareMathOperator{\rk}{rk}

\newtheorem{theorem}{Theorem}[section]
\newtheorem{prop}[theorem]{Proposition}
\newtheorem{lemma}[theorem]{Lemma}
\newtheorem{cor}[theorem]{Corollary}
\newtheorem{prob}[theorem]{Problem}
\theoremstyle{definition}

\theoremstyle{remark}
\newtheorem{rem}[theorem]{Remark}
\theoremstyle{plain}
\newtheorem*{theoremA}{Theorem A}
\newtheorem*{theoremB}{Theorem B}
\newtheorem*{theoremC}{Theorem C}
\newtheorem*{theoremD}{Theorem D}

\title[Freeness obstructions]{Freeness Obstructions for Self-Similar Groups Acting on Rooted Trees}

\author{Daniele D'Angeli}
\address{Dipartimento di Ingegneria, Universit\`a degli Studi Niccol\`o Cusano, Via Don Carlo Gnocchi 3, 00166 Roma, Italy}
\email{daniele.dangeli@unicusano.it}

\author{Emanuele Rodaro}
\address{Dipartimento di Matematica, Politecnico di Milano, Piazza Leonardo da Vinci 32, 20133 Milano, Italy}
\email{emanuele.rodaro@polimi.it}

\date{}

\begin{document}
\begin{abstract}
We prove two general obstructions to freeness for self-similar groups
and derive consequences for automaton groups. First, let $\mathcal A$
be a finite invertible Mealy automaton with a co-accessible identity
state. Then every finitely generated subgroup of $G(\mathcal A)$ that
is either self-similar and transitive on the first level, or
level-transitive on the rooted tree, is cyclic or non-free. This
answers a question of Grigorchuk and partially extends Sidki's
obstruction for automata of polynomial activity. The methods also
yield effective, quantitative, and geometric consequences for
relations and level Schreier graphs.

Second, we study the relation between freeness and bireversibility.
We prove that if a finite reduced reversible invertible Mealy
automaton generates a non-abelian free group, then its dual has a
bireversible connected component. If the states form a free basis and
the action is transitive on the first level, then the automaton itself
is bireversible. Finally, over a binary alphabet, the free-basis
assumption can also be removed: every finite reduced invertible Mealy
automaton generating a non-abelian free group and acting
level-transitively is bireversible.
\end{abstract}
\maketitle
\enlargethispage{4pt}
\section{Introduction}

Self-similar groups arise from actions on rooted trees that reproduce on every
rooted subtree. They include many of the standard examples in the theory of
groups of intermediate growth and of branch groups, and they are related to
profinite groups, symbolic dynamics, iterated monodromy groups, limit spaces,
and finite automata; see, for instance,
\cite{bhn:aut_til,fractalgr,dynamicssubgroup,rational,nekrashevyc}.
Despite the large number of examples and the variety of methods available for
the study of particular groups, rather little is known in general about the
algebraic structure of automaton groups. This is already visible in basic
decision problems, where positive results and undecidability phenomena occur
side by side; see \cite{DaRo14,freeness,israel,DaWeRo24} and the references therein.
Recent structural results for groups acting on rooted trees include
\cite{FrancoeurGarrido,LeBoudecMatteBon}.

A natural problem in this setting is to determine when a self-similar group can be free, or can contain a non-abelian free subgroup. Even for a group generated by a concrete finite automaton, proving freeness is usually specific to the example, while finding a single non-trivial relation may require a substantial computation. Aleshin constructed in 1983 a finite automaton whose states were later proved by M.~Vorobets and Y.~Vorobets to generate a free group of rank three \cite{Ale83,VV07}. Further constructions show that free groups of every finite rank can be generated by finite automata over a binary alphabet \cite{voro,SteVo}; affine and square-complex constructions provide related examples \cite{BruSid,Mozes,BoKi}, and free products of groups of order two occur as well \cite{Oli}.

There are also general mechanisms that prevent freeness. Nekrashevych proved that contracting self-similar groups contain no non-abelian free subgroup \cite{nekrashevycFree}, while Sidki established the same conclusion for groups generated by automata of polynomial activity \cite{sidki2}. The first part of this paper gives a further obstruction, based neither on contraction nor on an activity bound. It concerns automata possessing an identity state which is co-accessible, that is, reachable from every state in the transition graph. The following question was communicated to us by Grigorchuk \cite{privCom}.

\begin{prob}
Can a finite invertible Mealy automaton with a co-accessible identity
state generate a non-abelian free group whose action on the rooted tree
is level-transitive?
\end{prob}

We answer this question negatively in the following stronger form.

\begin{theoremA}
Let $\mathcal A$ be a finite invertible Mealy automaton over $A$.
Assume that a state $e$ induces the identity on $A^*$ and is reachable
from every state of $\mathcal A$. Let $H\leq G(\mathcal A)$ be finitely
generated. If either
\begin{enumerate}[label=(\alph*)]
\item $H$ is self-similar and transitive on $A$; or
\item $H$ is level-transitive on $A^*$,
\end{enumerate}
then $H$ is cyclic or non-free.
\end{theoremA}

The proof follows the sections of a finite generating set along random words. Co-accessibility forces these iterated sections to reach the identity with probability one. A counting argument then produces a non-trivial element in the kernel of a section homomorphism, and a free-group argument combines the local kernels arising on the different orbits into a defining relation. The construction is effective: finite graph computations and Stallings foldings give a practical obstruction to freeness and, when a rank drop occurs, an explicit relation. We also obtain an explicit, although coarse, upper bound on the length of a shortest relation. The same absorption mechanism has a geometric consequence: for a level-transitive action, asymptotic absorption of the generators into the identity forces the level Schreier graphs to have vanishing Cheeger constants. In particular, a level-transitive automaton group with a co-accessible identity state has no expanding level Schreier graphs and does not have Kazhdan's property~$(T)$; see Corollary~\ref{cor:coaccessible-no-T}.
\\
The second part of the paper concerns the role of bireversibility. Complete VH-square complexes connect bireversible automata with the theory of
lattices in products of trees; see \cite{Mozes} for the automata correspondence
and \cite{BurgerMozes,wise2} for the latter theory. This framework produces bireversible
automata with different algebraic behaviour, including
automaton groups that are free and automaton groups with Kazhdan's property~$(T)$ \cite{Mozes}. Thus bireversibility is compatible with a broad range of algebraic structures. Nevertheless, all known level-transitive constructions of non-abelian free automaton groups are bireversible. Starting from a bireversible example, one can easily obtain a non-bireversible presentation, but the resulting action need not remain transitive \cite{fragile}. This leads to the problem of determining to what extent transitivity and non-abelian freeness force bireversibility.

Our first result in this direction considers reversible automata. In this
case reversibility imposes strong constraints on the section maps, which
can be combined with the freeness obstruction developed above.

\begin{theoremB}
Let $\mathcal A$ be a finite reduced reversible invertible Mealy
automaton. If every connected component of the dual automaton
$\partial\mathcal A$ is non-bireversible, then $G(\mathcal A)$ is
cyclic or non-free.
\end{theoremB}

Equivalently, the dual of a finite reduced reversible automaton generating a non-abelian free group must contain a bireversible connected component. This necessary condition agrees with all known level-transitive constructions of free automaton groups. It also contrasts with the semigroup case, since reversible non-bireversible automata can generate free subsemigroups \cite{FraMi}.

When the states themselves form a free basis, the reversibility assumption can be removed. In this case a transition collision gives a fold in the Schreier graph of a first-level stabilizer, and the corresponding rank drop produces a non-trivial local kernel.

\begin{theoremC}
Let $\mathcal A=(X,A,\delta,\lambda)$ be a finite invertible Mealy
automaton with $|X|\geq2$. If $X$ is a free basis of
$G(\mathcal A)$ and $G(\mathcal A)$ acts transitively on $A$, then
$\mathcal A$ is bireversible.
\end{theoremC}

Thus, for first-level transitive automata, the remaining case is that of a reduced automaton whose states generate a free group with a redundant generating set. Over a binary alphabet, level-transitivity is sufficient to settle this case as well.
\begin{theoremD}
Let $\mathcal A$ be a finite reduced invertible Mealy automaton over a
binary alphabet. If $G(\mathcal A)$ is a non-abelian free group and
acts level-transitively, then $\mathcal A$ is bireversible.
\end{theoremD}

The arguments behind Theorems~A--D share a common local-to-global principle. For a finitely generated subgroup \(H\), kernels of section homomorphisms at suitable vertices are combined to produce a non-trivial relation among the generators. In the co-accessible case, these local kernels are obtained from absorption in a finite section graph; in the bireversibility results, they arise from section collisions or from a rank drop detected by Stallings folding.

Section~\ref{sec: preliminaries} fixes the notation and recalls the required facts on self-similar actions, free groups, and Mealy automata. Section~\ref{sec: section homomorphisms} develops the local-to-global criterion and its effective form. Sections~\ref{sec: probabilistic approach} and~\ref{sec: effective forms} prove Theorem~A and establish the effective and quantitative versions of the obstruction. Finally, Section~\ref{sec:reversible-not-bireversible} proves Theorems~B--D.

\section{Preliminaries}
\label{sec: preliminaries}

We recall the notation and the facts used below: the identities
for sections, basic properties of finitely generated free groups, and
the self-similar action defined by an invertible Mealy automaton.

\subsection{Actions on rooted trees}

Let $A$ be a finite non-empty alphabet.  The free monoid $A^*$ is the
rooted regular tree whose root is the empty word and whose edges join
$u$ to $ua$, with $u\in A^*$ and $a\in A$.  Its $n$-th level is
$A^n$.  An automorphism of $A^*$ fixes the root and preserves prefixes.
Throughout the paper actions are written on the right: the image of
$u$ under $g$ is denoted by $u\cdot g$, and
$u\cdot(gh)=(u\cdot g)\cdot h$.

For $g\in\Aut(A^*)$ and $u\in A^*$, there is a unique automorphism
$g|_u$, called the \emph{section of $g$ at $u$}, such that
\begin{equation}
 \label{eq: section definition}
 (uv)\cdot g=(u\cdot g)(v\cdot g|_u)
 \qquad(v\in A^*).
\end{equation}
The uniqueness follows by identifying the subtree rooted at $u$ with
$A^*$.  Repeated use of \eqref{eq: section definition} gives
\begin{equation}
 \label{eq: section identities}
 g|_{uv}=(g|_u)|_v,\qquad
 (gh)|_u=g|_u\,h|_{u\cdot g},\qquad
 g^{-1}|_u=\bigl(g|_{u\cdot g^{-1}}\bigr)^{-1}.
\end{equation}
Throughout the paper $G\leq\Aut(A^*)$ will denote a \emph{self-similar} group, that is, for every $g\in G$ and $u\in A^*$ we have $g|_u\in G$. We write $\Stab_G(u)$ for the
stabilizer of $u$ and $\Stab_G(A^n)$ for the pointwise stabilizer of
the $n$-th level.  The action is \emph{level-transitive} if it is
transitive on every $A^n$.

The action on the first level and the first-level sections determine an
automorphism.  More precisely, if $\tau_g\in\operatorname{Sym}(A)$ is
the permutation induced by $g$, then
$g\mapsto((g|_a)_{a\in A},\tau_g)$ is the usual embedding
$\Aut(A^*)\hookrightarrow\Aut(A^*)^A\rtimes\operatorname{Sym}(A)$.
The semidirect product is the permutational wreath product
$\Aut(A^*)\wr_A\operatorname{Sym}(A)$, where
$\operatorname{Sym}(A)$ permutes the coordinates indexed by $A$.
In particular, an automorphism fixing $A^n$ pointwise and having
trivial section at every vertex of $A^n$ is the identity.

Let $A=A_1\sqcup\cdots\sqcup A_k$ be the orbit decomposition of a
self-similar group $G$ on the first level.  Each $A_i^*$ is
$G$-invariant.  Indeed, if a word begins with a letter of $A_i$, then
its first letter remains in $A_i$ under the action of $G$; the action on
the remaining suffix is given by a section, which again belongs to
$G$.  Repeating the argument at each letter proves the claim.

\subsection{Free groups}

For a finite set $X$, let $F_X$ be the free group with basis $X$, and
put $\widetilde X=X\sqcup X^{-1}$.  Words are always understood as
freely reduced when they are viewed as elements of $F_X$.  If
$Y\subseteq F_X$, then $\langle Y\rangle$ denotes the subgroup
generated by $Y$; we call $Y$ a \emph{free basis} when the natural map
$F_Y\to\langle Y\rangle$ is an isomorphism.  The rank of a finitely
generated free group $L$ is denoted by $\rk(L)$. We use the convention
$[x,y]=xyx^{-1}y^{-1}$. For $x\in H$ and $L\leq H$, we write
\[
 C_H(x):=\{h\in H:hx=xh\},\qquad
 N_H(L):=\{h\in H:hLh^{-1}=L\},\qquad
 \Core_H(L):=\bigcap_{h\in H}hLh^{-1}
\]
for the centralizer of $x$, the normalizer of $L$, and the core of $L$,
respectively.

We use the following standard facts.  Free groups are torsion-free and
finitely generated free groups are Hopfian.  If $1\neq x$ belongs to a
free group, its centralizer is the unique maximal cyclic subgroup
containing $x$.  Non-trivial maximal cyclic subgroups are malnormal, root closed, and
self-normalizing; in particular, $N_F(C)=C$.  The
primitive root of a non-trivial word, and membership in the cyclic
subgroup that it generates, can be computed effectively.  If $L$ has finite index $d$ in a free group of finite rank
$q$, Schreier's formula gives
\begin{equation}
 \label{eq: Schreier formula}
 \rk(L)=1+d(q-1).
\end{equation}
We also use B.~H.~Neumann's coset-covering theorem: if a group is a
finite union of left cosets of subgroups, then at least one of those
subgroups has finite index \cite{NeumannCover}.

If a group generated by a finite set $W$ acts on the right on a finite
set $\Omega$, its symmetric Schreier graph has vertex set $\Omega$ and,
for every $s\in W\cup W^{-1}$, an edge
$\omega\xrightarrow{\,s\,}\omega\cdot s$.  The labels of the closed
paths based at $\omega$ form the stabilizer of $\omega$.  A spanning
tree, together with the non-tree edges, gives the usual Schreier basis
of this stabilizer.

For a finitely generated subgroup $L\leq F_X$, its Stallings automaton
is a finite inverse $X$-labelled graph recognizing the reduced words in
$L$ \cite{Sta}.  Folding computes this graph effectively.  In
particular, it computes a basis and the rank of $L$ and decides
membership.  We shall also use the following effective form.  Let $B$ be finite and
let $\varphi:F_B\to F_X$ be specified by the tuple
$(\varphi(b))_{b\in B}$.  Nielsen transformations reduce this tuple to
\[
 (c_1,\ldots,c_r,1,\ldots,1),
\]
where $\{c_1,\ldots,c_r\}$ is a basis of $\varphi(F_B)$.  Applying the
same transformations to the basis tuple $(b)_{b\in B}$ gives a basis
$(z_b)_{b\in B}$ of $F_B$.  If $r<|B|$, then one of the corresponding
$z_b$ satisfies $\varphi(z_b)=1$.  Equivalently,
this information can be recovered while folding the labelled bouquet
associated with the tuple $(\varphi(b))_{b\in B}$.

\subsection{Invertible Mealy automata}

A finite Mealy automaton is a tuple
$\mathcal A=(X,A,\delta,\lambda)$, where $X$ and $A$ are finite
non-empty sets, $X$ is the set of states, $A$ is the input and output
alphabet, and $\delta:X\times A\to X$ and
$\lambda:X\times A\to A$ are the transition and output maps.  We write
$x\xrightarrow{\,a\mid b\,}y$ when $\delta(x,a)=y$ and
$\lambda(x,a)=b$.  The automaton is \emph{invertible} if
$a\mapsto\lambda(x,a)$ is a permutation of $A$ for every state $x$.
In this case each state induces an automorphism of $A^*$ by
\begin{equation}
 \label{eq: automaton action}
 (av)\cdot x=(a\cdot x)(v\cdot x|_a),
 \qquad a\cdot x:=\lambda(x,a),
 \qquad x|_a:=\delta(x,a).
\end{equation}
The subgroup generated by the states is denoted by $G(\mathcal A)$
and is called the \emph{automaton group generated by $\mathcal A$}.
Equation~\eqref{eq: automaton action} shows that it is a finitely
generated self-similar subgroup of $\Aut(A^*)$. An automorphism
$g\in\Aut(A^*)$ is \emph{finite-state} if the set
$\{g|_u:u\in A^*\}$ is finite. Every state of a finite Mealy automaton
is finite-state.  The section identities show that products and inverses
of finite-state automorphisms are again finite-state; hence every element
of $G(\mathcal A)$ is finite-state. Throughout,
$\mathcal A$ denotes the chosen finite generating automaton, while
$G(\mathcal A)$ denotes the associated group. We refer to
\cite{nekrashevyc,rational} for general background. The automaton is \emph{reduced} if distinct states induce distinct
elements of $G(\mathcal A)$.  It is \emph{reversible} if $x\mapsto x|_a$ is a permutation of $X$ for every $a\in A$.  If
$\mathcal A$ is invertible, its inverse automaton has a state $x^{-1}$
for each $x\in X$ and a transition
$x^{-1}\xrightarrow{\,b\mid a\,}y^{-1}$ whenever
$x\xrightarrow{\,a\mid b\,}y$ is a transition of $\mathcal A$.  The
automaton is \emph{bireversible} if both $\mathcal A$ and its inverse
are reversible. A state $e$ is an \emph{identity state} if it induces the identity of
$A^*$.  It is \emph{co-accessible} if, from every state $x$, some input
word labels a path from $x$ to $e$.  Equivalently, every state has a
section equal to $e$. The \emph{dual automaton} $\partial\mathcal A$ has state set $A$ and
alphabet $X$, with a transition
$a\xrightarrow{\,x\mid y\,}b$ whenever
$x\xrightarrow{\,a\mid b\,}y$ is a transition of $\mathcal A$.  If $\mathcal A$ is invertible, each input letter $x$ acts as a
permutation of the state set $A$ of the dual. Hence the connected
components of $\partial\mathcal A$ are precisely the orbits of
$G(\mathcal A)$ on $A$. For a component $C\subseteq A$, the
transitions whose source vertex lies in $C$ form a reversible Mealy
automaton with state set $C$ and alphabet $X$, not necessarily an
invertible one.
We call the component \emph{bireversible} when this induced Mealy
automaton is invertible and bireversible. 

\section{Section homomorphisms and free presentations}
\label{sec: section homomorphisms}

The purpose of this section is to convert kernels of section maps into
relations among a chosen generating set.  Let
$G\leq\operatorname{Aut}(A^*)$ be self-similar, let $H\leq G$ be
finitely generated, and let $W$ be a finite generating set of $H$.
Throughout this section,
\[
 \rho:F_W\to H
\]
denotes the natural epimorphism. For $u\in A^*$, put
\[
 P_u:=\rho^{-1}(\operatorname{Stab}_H(u)),\qquad
 \sigma_u:P_u\to G,\quad \sigma_u(w):=\rho(w)|_u.
\]
The map \(\sigma_u\) is called the \emph{section homomorphism} at
\(u\). It is a standard fact in the theory of self-similar groups (see,
e.g., \cite{nekrashevyc}) that the restriction of the section map
to a vertex stabilizer is a group homomorphism. Since we shall use this
property repeatedly, and because our formulation takes values in the
ambient group \(G\) rather than in a self-similar subgroup, we include
the short proof for completeness.

\begin{lemma}
\label{lem:section-homomorphism}
For every \(u\in A^*\), the map \(\sigma_u:P_u\longrightarrow G\) is a homomorphism.
\end{lemma}

\begin{proof}
Let \(w,z\in P_u\). Since both \(\rho(w)\) and \(\rho(z)\) fix \(u\),
also \(wz\in P_u\). Using the product identity 
\eqref{eq: section identities} we get
\[
\rho(wz)|_u
=
\rho(w)|_u\,\rho(z)|_{u\cdot\rho(w)}
=
\rho(w)|_u\,\rho(z)|_u.
\]
Hence \(\sigma_u(wz)=\sigma_u(w)\sigma_u(z),\)
so \(\sigma_u\) is a homomorphism.
\end{proof}
For $n\geq1$ put
\[
 P_n:=\rho^{-1}(\operatorname{Stab}_H(A^n))
     =\bigcap_{u\in A^n}P_u,
 \qquad
 \Sigma_n(w):=(\sigma_u(w))_{u\in A^n}.
\]
Thus $\Sigma_n:P_n\to G^{A^n}$ is the product of the section
homomorphisms over the $n$-th level.  The next proposition identifies
its kernel.

\begin{prop}
\label{prop: level-kernel}
For every $n\geq1$,
\[
 \ker(\Sigma_n)=\ker(\rho).
\]
In particular, $W$ is not a free basis of $H$ if and only if
$\Sigma_n$ is not injective for some, equivalently every, $n\geq1$.
\end{prop}

\begin{proof}
If $w\in\ker(\rho)$, then $\rho(w)$ fixes every vertex and all its
sections are trivial; hence $w\in\ker(\Sigma_n)$.

Conversely, let $w\in\ker(\Sigma_n)$.  By definition, $w\in P_n$, so
$\rho(w)$ fixes every vertex $u\in A^n$.  Moreover,
$\sigma_u(w)=\rho(w)|_u=1$ for every such $u$.  Every vertex of level at
least $n$ has a unique form $uv$ with $u\in A^n$ and $v\in A^*$, and
\eqref{eq: section definition} gives
\[
 (uv)\cdot\rho(w)
 =(u\cdot\rho(w))\bigl(v\cdot\rho(w)|_u\bigr)=uv.
\]
Vertices of levels smaller than $n$ are prefixes of vertices in $A^n$;
since tree automorphisms preserve prefixes, they are fixed as well.
Thus $\rho(w)$ fixes all of $A^*$ and is the identity.
\end{proof}

\subsection{Patching local kernels}

The methods used later produce, for each orbit of a level, a
non-trivial element in the kernel of a section homomorphism at one
vertex of that orbit.  These elements need not coincide.  The next two
lemmas combine them into a single element of $\ker(\Sigma_n)$.  The
first keeps the successive commutators non-trivial; the second performs
the patching.

\begin{lemma}
\label{lem: non trivial commutators}
Let $F_B$ be the free group on a finite basis $B$ with $|B|\geq2$,
and let $x,y\in F_B\setminus\{1\}$. There exists
$v\in\{1\}\cup B$ such that
\[
 [v^m x v^{-m},y]\neq1
 \qquad
 \text{for every }m\in\mathbb Z\setminus\{0\}.
\]
In particular, $v$ can be found effectively.
\end{lemma}

\begin{proof}
If $[x,y]\neq1$, take $v=1$.  Assume that $[x,y]=1$.  Then
\[
 C:=C_{F_B}(x)=C_{F_B}(y)
\]
is the maximal cyclic subgroup containing $x$ and $y$.  Since $F_B$ is
not cyclic, some $v\in B$ lies outside $C$. Suppose that $[v^m x v^{-m},y]=1$ for some $m\neq0$.  Then
\[
 C_{F_B}(v^m x v^{-m})=v^mCv^{-m}=C_{F_B}(y)=C.
\]
Thus $v^m\in N_{F_B}(C)=C$.  Since $C$ is root closed, $v\in C$, a
contradiction. For the effective assertion, compute the primitive root $r$ of $x$, so
that $C=\langle r\rangle$, and test the elements of the finite basis
$B$ for membership in $C$.  The first element outside $C$ is an
admissible choice of $v$.
\end{proof}

% We retain the notation $G,H,W,\rho,P_u,\sigma_u$ and $\Sigma_n$ fixed
% at the beginning of the section.

\begin{lemma}
\label{lem: local-kernel-patching}
Assume $\operatorname{rk}(F_W)\geq2$, fix $n\geq1$, and let
$A^n=O_1\sqcup\cdots\sqcup O_s$ be the decomposition into
$H$-orbits. The following are equivalent:
\begin{enumerate}[label=(\roman*)]
 \item for every $j$ there exists $u_j\in O_j$ such that
       $\ker(\sigma_{u_j})\neq1$;
 \item $\ker(\Sigma_n)\neq1$.
\end{enumerate}
Given the finite action on $A^n$ and explicit local kernel words, the
construction of an element of $\ker(\Sigma_n)\setminus\{1\}$ is
effective.
\end{lemma}

\begin{proof}
Assume (ii), and choose
$1\neq q\in\ker(\Sigma_n)$.  Then
$q\in\ker(\sigma_u)$ for every $u\in A^n$, so (i) follows.  We
now assume (i). The proof has three steps: propagate one local kernel along each orbit,
move all resulting words into the stabilizer of the whole level, and
then eliminate the section coordinates one at a time.  We first
propagate a local kernel along an orbit.  Let $u,v$ lie in the
same $H$-orbit and choose $t\in F_W$ with
$v\cdot\rho(t)=u$.  If $1\neq z\in\ker(\sigma_u)$, then
$z_v:=tzt^{-1}$ is non-trivial.  Moreover, $z_v$ fixes $v$, because
\[
 v\cdot\rho(tzt^{-1})
 =u\cdot\rho(zt^{-1})=u\cdot\rho(t^{-1})=v.
\]
Using the product and inverse identities for sections, and the facts
that $\rho(z)$ fixes $u$ and $\rho(z)|_u=1$, we obtain
\[
 \rho(z_v)|_v
 =\rho(t)|_v\,\rho(z)|_u\,\rho(t^{-1})|_u=1.
\]
Thus $1\neq z_v\in\ker(\sigma_v)$.  We may therefore choose, for every
$u\in A^n$, an element $1\neq z_u\in\ker(\sigma_u)$.

We next place all local kernel elements in the common domain $P_n$.
Let $N$ be the exponent of the finite permutation group induced by $H$
on $A^n$.  For every $q\in F_W$, the element $\rho(q)^N$ fixes $A^n$
pointwise, so $q^N\in P_n$.  Replacing $z_u$ by $z_u^N$ preserves
non-triviality, since $F_W$ is torsion-free, and preserves membership in
$\ker(\sigma_u)$, since $\sigma_u$ is a homomorphism.  Hence we may
assume
\[
 1\neq z_u\in P_n\cap\ker(\sigma_u)
 \qquad(u\in A^n).
\]

It remains to kill the coordinates one at a time without producing the
trivial word.  Enumerate $A^n=\{u_1,\ldots,u_r\}$.  We construct
inductively a non-trivial $q_j\in P_n$ such that
$\sigma_{u_i}(q_j)=1$ for $1\leq i\leq j$.  Set $q_1=z_{u_1}$.
Suppose $q_j$ is defined.  By
Lemma~\ref{lem: non trivial commutators}, choose $v\in F_W$ so that
\[
 q_{j+1}:=[v^N z_{u_{j+1}}v^{-N},q_j]\neq1.
\]
The four entries in this commutator lie in $P_n$, hence so does
$q_{j+1}$.  The commutator is arranged so that one of its entries is
trivial under each section homomorphism already treated.  If $i\leq j$,
the second entry $q_j$ has trivial $\sigma_{u_i}$-image, so
$\sigma_{u_i}(q_{j+1})=1$.  At $u_{j+1}$, the first entry has trivial
image because
\[
 \sigma_{u_{j+1}}(v^N z_{u_{j+1}}v^{-N})
 =\sigma_{u_{j+1}}(v^N)\,1\,
  \sigma_{u_{j+1}}(v^N)^{-1}=1.
\]
Thus $q_{j+1}$ has the required property.  At the end,
$1\neq q_r\in\ker(\Sigma_n)$.

All steps are effective.  The finite action on $A^n$ gives $N$, orbit
representatives and transport words can be found in the finite
Schreier graph, and Lemma~\ref{lem: non trivial commutators} gives an
effective choice of the conjugating elements.
\end{proof}

\begin{theorem}
\label{theo: local-global-criterion}
Assume $\operatorname{rk}(F_W)\geq2$ and fix $n\geq1$. The following are
equivalent:
\begin{enumerate}[label=(\roman*)]
 \item $W$ is not a free basis of $H$;
 \item for every $H$-orbit $O\subseteq A^n$ there exists $u\in O$ such
       that $\ker(\sigma_u)\neq1$;
 \item $\ker(\Sigma_n)\neq1$.
\end{enumerate}
\end{theorem}

\begin{proof}
By Proposition~\ref{prop: level-kernel},
$W$ fails to be a free basis exactly when the global map $\Sigma_n$ has
a non-trivial kernel; this proves (i)$\Leftrightarrow$(iii).  By
Lemma~\ref{lem: local-kernel-patching}, the latter happens exactly when
one local kernel is present on each orbit of the $n$-th level, proving
(ii)$\Leftrightarrow$(iii).
\end{proof}
Thus it remains to find, on each orbit, one vertex $u$ with
$\ker(\sigma_u)\neq1$; the patching lemma then produces a non-trivial
element of $\ker(\rho)$.

\subsection{Self-similar subgroups}

Assume now that $H$ is self-similar. Then $\sigma_u(P_u)\leq H$ for every
$u\in A^*$.  This additional invariance allows us to move a local
kernel from an arbitrary depth to the first level: one follows the
successive sections of a non-trivial element until the first time they
become trivial.

\begin{lemma}
\label{lem: depth-reduction}
If $\sigma_u$ is not injective for some
$u=a_1\cdots a_m\in A^*$, then $\sigma_{a_j}$ is not injective for some
$j\in\{1,\ldots,m\}$.
\end{lemma}

\begin{proof}
Choose $1\neq w\in\ker(\sigma_u)$ and set $g_0:=\rho(w)$.  Since
$w\in P_u$, the element $g_0$ fixes the word $u$ and hence every prefix
of $u$.  If $g_0=1$, then $w$ belongs to $P_{a_1}$ and
$\sigma_{a_1}(w)=1$, so there is nothing to prove.

Assume $g_0\neq1$.  Define successively
$g_j:=g_{j-1}|_{a_j}$ for $1\leq j\leq m$.  Because $g_0$ fixes
$a_1\cdots a_m$, the element $g_{j-1}$ fixes $a_j$ and its section
$g_j$ fixes the remaining suffix.  Self-similarity gives $g_j\in H$,
and $g_m=g_0|_u=1$.  Let $j$ be minimal with $g_j=1$.  Then
$g_{j-1}\neq1$.  Choose $z\in F_W$ with $\rho(z)=g_{j-1}$; necessarily
$z\neq1$.  Since $g_{j-1}$ fixes $a_j$, we have $z\in P_{a_j}$, and
$\sigma_{a_j}(z)=g_{j-1}|_{a_j}=g_j=1$.  Thus $\sigma_{a_j}$ is not
injective.
\end{proof}

Let $A=A_1\sqcup\cdots\sqcup A_k$ be the decomposition into
$H$-orbits. Since $H$ is self-similar, each
$A_i^*$ is $H$-invariant.

\begin{cor}
\label{cor: first level reduction}
Assume $\operatorname{rk}(F_W)\geq2$. The following are equivalent:
\begin{enumerate}[label=(\roman*)]
 \item $W$ is not a free basis of $H$;
 \item for every $i$ there exist $m_i\geq1$ and $u_i\in A_i^{m_i}$ such
       that $\ker(\sigma_{u_i})\neq1$;
 \item for every $i$ there exists $a_i\in A_i$ such that
       $\ker(\sigma_{a_i})\neq1$.
\end{enumerate}
\end{cor}

\begin{proof}
Assume (i) and choose $1\neq w\in\ker(\rho)$.  For every $i$, choose
any $u_i\in A_i$.  The element $w$ fixes $u_i$ and has trivial section
there, so $1\neq w\in\ker(\sigma_{u_i})$; this proves (ii), with
$m_i=1$.  The implication (ii)$\Rightarrow$(iii) is
Lemma~\ref{lem: depth-reduction}.  Finally, the vertices
$a_i$ in (iii) meet every $H$-orbit on $A$, so
Theorem~\ref{theo: local-global-criterion} with $n=1$ gives a
non-trivial element of $\ker(\rho)$ and proves (i).
\end{proof}

\subsection{Lifted section homomorphisms}

Assume that $H$ is self-similar.  Then each $\sigma_a$ has codomain
$H$, but deciding its kernel requires knowledge of the relations in
$H$, which are precisely what we seek.  To obtain a condition that can
be checked inside the free group $F_W$, we choose words representing
the sections of the generators and lift $\sigma_a$ to a homomorphism
$P_a\to F_W$.  A kernel of the lifted map is a kernel of $\sigma_a$;
its rank can be computed by Stallings folding.

For every $a\in A$ and $x\in W$,
choose a word $\lambda(a,x)\in F_W$ representing the section
$\rho(x)|_a$.  Define the representatives for inverse letters by
\[
 \lambda(a,x^{-1})
 :=\lambda(a\cdot\rho(x)^{-1},x)^{-1}.
\]
This is compatible with the inverse identity in
\eqref{eq: section identities}. For a word $z=x_1\cdots x_\ell\in\widetilde W^*$, follow the orbit of
$a$ while the letters of $z$ are read: put $a_0=a$ and
$a_j=a_{j-1}\cdot\rho(x_j)$.  The lifted output is the product of the
chosen representatives encountered along this path:
\[
 \widehat\sigma_a(z)
 :=\lambda(a_0,x_1)\lambda(a_1,x_2)\cdots
   \lambda(a_{\ell-1},x_\ell).
\]
\begin{lemma}
\label{lem: lifted section map}
The formula above defines a well-defined map
$\widehat\sigma_a:F_W\to F_W$.  For all $z,t\in F_W$,
\[
 \rho(\widehat\sigma_a(z))=\rho(z)|_a,
 \qquad
 \widehat\sigma_a(zt)
 =\widehat\sigma_a(z)\widehat\sigma_{a\cdot\rho(z)}(t).
\]
Consequently, the restriction
$\widehat\sigma_a:P_a\to F_W$ is a homomorphism and
$\rho\circ\widehat\sigma_a=\sigma_a$.
\end{lemma}

\begin{proof}
Suppose that, while reading a word from the state $b\in A$, two
consecutive letters are $x,x^{-1}$ with $x\in W$.  Their contribution
to the lifted output is
\[
 \lambda(b,x)\lambda(b\cdot\rho(x),x^{-1})
 =\lambda(b,x)\lambda(b,x)^{-1}=1.
\]
The same calculation applies to $x^{-1}x$.  Hence insertion or deletion
of a freely cancelling pair does not change
$\widehat\sigma_a(z)$, and the construction descends from words to
$F_W$.

The first identity follows by induction on the length of $z$.  For one
letter it is the definition of $\lambda$, including the inverse case.
If $z=z'x$, the induction hypothesis and the product identity for
sections give
\[
 \rho(z)|_a
 =\rho(z')|_a\,\rho(x)|_{a\cdot\rho(z')}
 =\rho\bigl(\widehat\sigma_a(z')
   \lambda(a\cdot\rho(z'),x)\bigr).
\]
This is the required formula.  Splitting the defining product after the letters of $z$ gives
\[
 \widehat\sigma_a(zt)
 =\widehat\sigma_a(z)\widehat\sigma_{a\cdot\rho(z)}(t).
\]
If $z\in P_a$, then $a\cdot\rho(z)=a$, so the preceding identity
becomes
$\widehat\sigma_a(zt)=\widehat\sigma_a(z)\widehat\sigma_a(t)$ for
$z,t\in P_a$.  The first identity then gives
$\rho\circ\widehat\sigma_a=\sigma_a$ on $P_a$.
\end{proof}

\begin{prop}
\label{prop: lifted rank criterion}
Let $A=A_1\sqcup\cdots\sqcup A_k$ be the decomposition into
$H$-orbits, and choose $a_i\in A_i$.  For each $i$, the restriction
$\widehat\sigma_{a_i}:P_{a_i}\to F_W$ is injective if and only if
\[
 \rk(\widehat\sigma_{a_i}(P_{a_i}))=\rk(P_{a_i}).
\]
In particular, if the inequality is strict for every $i$, then $W$ is
not a free basis of $H$.  Given the finite action on $A$ and the chosen
section representatives $\lambda(a,x)$, these rank conditions are
decidable; when they are strict on every orbit, a non-trivial word in
$\ker(\rho)$ can be constructed.
\end{prop}

\begin{proof}
The subgroup $P_{a_i}$ has finite index in $F_W$, because it is the
stabilizer of a point in a finite $F_W$-set.  It is therefore a
finitely generated free group.  The map $\widehat\sigma_{a_i}:P_{a_i}\to
\widehat\sigma_{a_i}(P_{a_i})$ is onto.  If it is injective, the two
free groups are isomorphic and have the same rank.  Conversely, suppose
that their ranks are equal.  Choose an isomorphism from
$\widehat\sigma_{a_i}(P_{a_i})$ back to $P_{a_i}$ and compose it with
$\widehat\sigma_{a_i}$.  The resulting map is a surjective
endomorphism of the finitely generated free group $P_{a_i}$, hence is
injective by Hopficity.  Therefore $\widehat\sigma_{a_i}$ itself is
injective.  Thus strict rank inequality is equivalent to a non-trivial
kernel.  Such a kernel is contained in $\ker(\sigma_{a_i})$ by
Lemma~\ref{lem: lifted section map}.  If the inequality is strict on
every orbit, Corollary~\ref{cor: first level reduction} shows that
$W$ is not a free basis.

This also proves the last assertion: the finite action computes the
subgroups $P_{a_i}$, the chosen representatives $\lambda(a,x)$ compute
the lifted images of their generators, and Stallings folding computes
the ranks of these images.  The construction of a relation is described
below.
\end{proof}

\subsection{Computing a relation from a rank drop}

The preceding proposition gives a finite sufficient test for non-freeness.  A rank drop
on every first-level orbit yields explicit elements of the relevant
section kernels, which can then be patched into a defining relation.
Assume that the action of $F_W$ on $A$ and the words $\lambda(a,x)$
are given.  For each orbit $A_i$, proceed as follows.

\begin{enumerate}[label=(\arabic*),leftmargin=2.2em]
\item Choose $a_i\in A_i$ and construct the Schreier graph of the action
      on $A_i$, rooted at $a_i$.  A spanning tree gives a Schreier basis
      $B_i$ of $P_{a_i}$.
\item For every $b\in B_i$, compute $\widehat\sigma_{a_i}(b)$ by following
      the corresponding path in the Schreier graph and replacing each
      letter by the chosen representative of its section.  This gives
      the image of a basis of $P_{a_i}$ under the lifted homomorphism.
\item Form a based bouquet with one loop for each $b\in B_i$, subdivide
      that loop according to the reduced word
      $\widehat\sigma_{a_i}(b)$, and fold the resulting $W$-labelled
      graph.  Its core is the Stallings graph of
      $\widehat\sigma_{a_i}(P_{a_i})$, hence determines its rank.
\item Apply Nielsen reduction to the tuple
      $(\widehat\sigma_{a_i}(b))_{b\in B_i}$ and perform the same
      elementary Nielsen moves on the basis tuple $(b)_{b\in B_i}$.
      If the image rank is smaller than $|B_i|$, the reduced image tuple
      contains an entry equal to $1$; the corresponding entry of the
      transformed basis tuple is a non-trivial word in
      $\ker(\widehat\sigma_{a_i})$.  Substituting the Schreier basis
      words gives $1\neq z_i\in P_{a_i}\leq F_W$.  If this happens for
      every orbit, Lemma~\ref{lem: local-kernel-patching} applied to the
      words $z_i$ produces a non-trivial element of $\ker(\rho)$.
\end{enumerate}

Every step is a finite computation.  When the rank drops on every
orbit, the procedure returns an explicit defining relation.  If the
rank does not drop on some orbit, the test gives no conclusion, since
$\ker(\sigma_{a_i})$ may be non-trivial even when
$\ker(\widehat\sigma_{a_i})=1$. Figure~\ref{fig:rank-test} illustrates the procedure on a two-state
binary automaton.  It also shows why the lifted maps are retained in the
effective part of the argument: their images can be computed by Stallings folding, and a rank drop can be converted into an explicit
relation.

\begin{figure}[ht]
\centering
\resizebox{\textwidth}{!}{%
\begin{tikzpicture}[
  >=Stealth,
  state/.style={circle,draw,minimum size=8mm,inner sep=1pt},
  every edge/.style={draw,->},
  every loop/.style={min distance=8mm}
]
% The automaton and its dual
\node[state] (qa) at (0,2.0) {$a$};
\node[state] (qb) at (3.0,2.0) {$b$};
\node at (-0.15,2.85) {$\mathcal A$};
\path
 (qa) edge[bend left=19] node[above] {$\substack{0\mid0\\[-1pt]1\mid1}$} (qb)
 (qb) edge[bend left=19] node[below] {$0\mid1$} (qa)
 (qb) edge[loop right] node[right] {$1\mid0$} (qb);

\node[state] (q0) at (0,-1.1) {$0$};
\node[state] (q1) at (3.0,-1.1) {$1$};
\node at (0.15,-2.05) {$\partial\mathcal A$};
\path
 (q0) edge[loop left] node[left] {$a\mid b$} (q0)
 (q1) edge[loop right] node[right] {$a\mid b$} (q1)
 (q0) edge[bend right=24] node[below] {$b\mid a$} (q1)
 (q1) edge[bend right=24] node[above] {$b\mid b$} (q0);

% Stallings automaton of P_0
\node at (7,3.0) {$\mathcal S(P_0)$};
\node[state] (p) at (6.2,1.8) {$p$};
\node[state] (q) at (6.2,-0.6) {$q$};
\path
 (p) edge[loop above] node[above] {$a$} (p)
 (q) edge[loop below] node[below] {$a$} (q)
 (p) edge[bend left=17] node[right] {$b$} (q)
 (q) edge[bend left=17] node[left] {$b$} (p);

% Stallings automaton of the image
\node at (9.55,2.45) {$\mathcal S(\widehat\sigma_0(P_0))$};
\node[state] (r) at (9.55,0.65) {$\ast$};
\path
 (r) edge[loop above] node[above] {$a$} (r)
 (r) edge[loop right] node[right] {$b$} (r);
\end{tikzpicture}%
}
\caption{An effective rank test. On the left are the automaton and its
dual.  On the right are the Stallings automata of the first-level
stabilizer preimage $P_0$ and of its image under the lifted section
map.  Their ranks are $3$ and $2$, respectively.}
\label{fig:rank-test}
\end{figure}
For Figure~\ref{fig:rank-test}, \(P_0\) has Schreier basis
\(x=a\), \(y=bab^{-1}\), \(z=b^2\), with
\[
\widehat{\sigma}_0(x)=b,\qquad
\widehat{\sigma}_0(y)=aba^{-1},\qquad
\widehat{\sigma}_0(z)=ab.
\]
Hence \(\operatorname{rk}(P_0)=3\) and
\(\operatorname{rk}(\widehat{\sigma}_0(P_0))=2\), while
\(yzx^{-1}z^{-1}=bab\,a^{-1}b^{-2}\in\ker(\widehat{\sigma}_0)\).
The local-to-global construction yields, for instance,
\[
[a^{-1}ba^{-1}b^{-1},(ba^{-1})^2b^{-2}]=1.
\]

\section{A probabilistic criterion for non-freeness}
\label{sec: probabilistic approach}
We next use the dynamics of sections to produce the local kernels needed in the previous
section. On a fixed level, we count the section-trivial edges.If their number is at least the number of
vertices, they contain a reduced cycle whose label is a non-trivial
word with trivial section.

Let \(H=\langle W\rangle\leq G\le \Aut(A^*)\) be a finitely generated group, where \(W\subseteq G\) is a finite generating set, and let
\[
\pi:F_W\twoheadrightarrow H
\]
be the natural epimorphism. Fix \(m\ge1\). For \(w\in W\) and \(u\in A^m\), consider the edge
\(u\xrightarrow{\,w\,}u\cdot\pi(w)\) in the Schreier graph of the action
of \(F_W\) on \(A^m\). It is \emph{section-trivial} if
\(\pi(w)|_u=1\). Let \(T_w(m)\) and \(E_w(m)\) denote the sets of
section-trivial and non-section-trivial edges labelled by \(w\), and put
\[
T(m):=\bigsqcup_{w\in W}T_w(m),
\qquad
E(m):=\bigsqcup_{w\in W}E_w(m).
\]
Then
\[
|T(m)|+|E(m)|=|W|\,|A|^m.
\]
\begin{lemma}
\label{lem:section-trivial-cycle}
If \( |T(m)|\geq |A|^m\), then there exist \(u\in A^m\) and
\(1\neq z\in F_W\) such that
\[
u\cdot\pi(z)=u,
\qquad
\pi(z)|_u=1.
\]
\end{lemma}

\begin{proof}
Let \(\Gamma_m\) be the undirected multigraph obtained from the Schreier
graph of the action of \(F_W\) on \(A^m\) by retaining only the edges
in \(T(m)\) and forgetting their orientations. Loops and parallel edges
are retained. Since \(|E(\Gamma_m)|\geq |V(\Gamma_m)|,\) some connected component contains a cyclically reduced closed path
\(c\), based at \(u\). Label each traversal of an edge of \(c\) by \(w\) or \(w^{-1}\),
according to whether it agrees with or opposes the corresponding edge
of \(T(m)\), and let \(z\in\widetilde W^*\) be the resulting word.
Since each \(\pi(w)\) permutes \(A^m\), two consecutive inverse
letters, including the last and the first, would correspond to an
immediate reversal of the same edge. Thus \(z\) is freely reduced and
\(z\neq1\) in \(F_W\). Since \(c\) is closed,
\[
u\cdot\pi(z)=u.
\]
Every forward traversal has trivial section at its initial vertex. For
a reversed edge \(p\xrightarrow{w}q\), this follows from \(\pi(w)^{-1}|_q=(\pi(w)|_p)^{-1}=1.\) The product identity for sections now gives \(\pi(z)|_u=1.\)
\end{proof}
The \emph{section graph} \(\mathcal S_H(W)\) is the deterministic
\(A\)-labelled graph generated by \(\pi(W)\) under sections. Its vertex
set is
\[
V(\mathcal S_H(W))
=
\{\pi(w)|_u:w\in W,\ u\in A^*\},
\]
and, for every \(g\in V(\mathcal S_H(W))\) and \(a\in A\), it has the edge \(g\xrightarrow{\,a\,}g|_a.\) Note that for a fixed $w\in W$, we have
$|T_w(m)|=|\{u\in A^m:\pi(w)|_u=1\}|$. Thus we may define
\[
 \chi_W(m)
 :=\frac{|T(m)|}{|W|\,|A|^m}
 =\frac{1}{|W|}\sum_{w\in W}
   \frac{|\{u\in A^m:\pi(w)|_u=1\}|}{|A|^m}
\]
Equivalently, \(\chi_W(m)\) is the probability that a uniformly chosen path of length \(m\) in \(A^*\), starting at a uniformly chosen
generator of $W$, ends at the identity. Since \(1|_a=1\) for every \(a\in A\),
the sequence \(\chi_W(m)\) is non-decreasing, so the limit \(\chi_W=\lim_{m\to\infty}\chi_W(m)\) exists, and this may be interpreted as the probability that a random walk starting at a vertex of $W$ eventually ends in the identity state. In case $H$ is self-similar and $A=A_1\sqcup\cdots\sqcup A_k$ is its first-level orbit decomposition, then $H$ acts on the sub-tree $A_i^*$, in this case we put an apex to the section graph, so  \(\mathcal S_H^{(i)}(W)\) is the deterministic
\(A_i\)-labelled graph generated by \(\pi(W)\) under sections taken on the alphabet $A_i$. So its vertex set is
\(
V(\mathcal S_H^{(i)}(W))
=
\{\pi(w)|_u:w\in W,\ u\in A_i^*\},
\)
and, for every \(g\in V(\mathcal S_H^{(i)}(W))\) and \(a\in A_i\) , it has the edge \(g\xrightarrow{\,a\,}g|_a.\) Similarly as before we have
$|T^{(i)}_w(m)|=|\{u\in A_i^m:\pi(w)|_u=1\}|$, \(T^{(i)}(m):=\bigsqcup_{w\in W}T^{(i)}_w(m)\), and the probability of eventually being absorbed to the identity using elements in the alphabet $A_i$ is defined by
\[
 \chi^{(i)}_W(m)
 :=\frac{|T^{(i)}(m)|}{|W|\,|A_i|^m}
 =\frac{1}{|W|}\sum_{w\in W}
   \frac{|\{u\in A_i^m:\pi(w)|_u=1\}|}{|A_i|^m}
\]
also in this case the limit
\(\chi_W^{(i)}=\lim_{m\to\infty}\chi_W^{(i)}(m)\) exists. Henceforth, when the apex $i$ is present it is implicitly stated that we are considering a self-similar group $H$. 
\begin{prop}\label{pro:probabilistic-criterion}
  Suppose that $H=\langle W\rangle$, $|
  W|\ge 2$ is self-similar and let $A=A_1\sqcup\cdots\sqcup A_k$ be its first-level orbit decomposition.   If
\(
\chi_W^{(i)}>\frac1{|W|}
\)
for $i\in [1,k]$ then the natural epimorphism
\(\pi:F_W\to H\) is not injective.
\end{prop}

\begin{proof}
For each \(i\), choose \(m_i\) such that
\(\chi_W^{(i)}(m_i)>1/|W|\). Then
\(|T^{(i)}(m_i)|>|A_i|^{m_i}\), so
Lemma~\ref{lem:section-trivial-cycle} gives
\(u_i\in A_i^{m_i}\) and \(1\neq z_i\in F_W\) with
\(u_i\cdot\pi(z_i)=u_i\) and \(\pi(z_i)|_{u_i}=1\).
Hence \(1\neq z_i\in\ker(\sigma_{u_i})\) for every first-level orbit.
Corollary~\ref{cor: first level reduction} yields
\(\ker(\pi)\neq1\).
\end{proof}
The threshold \(1/|W|\) is exactly the value at which the number of
section-trivial edges exceeds the number of vertices and hence forces a
cycle. The probabilistic criterion becomes especially simple when the section
graph is finite.  If the identity is reachable from every section,
then a random sequence of letters reaches it with probability one. The following is a standard fact in the theory of Markov chains.

\begin{lemma}
\label{lem: absorption-one}
With the above notation. Assume that
$\mathcal S_H(W)$ is finite and that $1$ is reachable from every vertex.  Then a random walk whose labels are chosen independently and
uniformly in $A$ reaches $1$ with probability one, that is, \(\chi_W=1\).
\end{lemma}

\begin{proof}
The identity is absorbing because $1|_a=1$ for every $a\in A$. Since the graph is finite and every vertex reaches $1$, there is an
integer $L$ such that from each vertex one can reach $1$ by a word of
length at most $L$.  Let $d:=|A|$.  Starting from any vertex, the probability of following one fixed such word is at least $d^{-L}$;
if the word has length smaller than $L$, the remaining steps may be
arbitrary because $1$ is absorbing.  Consequently, the probability of avoiding $1$
for $rL$ steps is at most $(1-d^{-L})^r$, which tends to zero as
$r\to\infty$.  Hence the walk reaches $1$ with probability one.  Averaging over the
finite set of initial states $\pi(W)$ gives $\chi_W=1$.
\end{proof}
Finiteness and reachability is a property preserved by taking finitely generated subgroups. 
\begin{lemma}
\label{lem:finite-section-inheritance}
Let \(Y\subseteq\langle W\rangle=H\leq\operatorname{Aut}(A^*)\) be finite.
\begin{enumerate}[label=(\roman*)]
 \item If \(\mathcal S_H(W)\) is finite, then
       \(\mathcal S_H(Y)\) is finite.
 \item If \(1\) is reachable from every vertex of
       \(\mathcal S_H(W)\), then the same holds for
       \(\mathcal S_H(Y)\).
\end{enumerate}
If \(H\) is self-similar and \(A_i\) is a first-level \(H\)-orbit,
the same assertions hold for
\(\mathcal S_H^{(i)}(W)\) and \(\mathcal S_H^{(i)}(Y)\).
\end{lemma}

\begin{proof}
Set
\[
C:=\{w|_u:w\in W,\ u\in A^*\}
   =V\bigl(\mathcal S_H(W)\bigr).
\]
For \(w\in W\) and \(u\in A^*\),
\[
(w^{-1})|_u
=
\bigl(w|_{u\cdot w^{-1}}\bigr)^{-1}
\in C^{-1}.
\]
Let \(y\in Y\), and write
\[
y=w_1^{\varepsilon_1}\cdots w_\ell^{\varepsilon_\ell},
\qquad
w_j\in W,\quad \varepsilon_j\in\{\pm1\}.
\]
For every \(u\in A^*\), the product identity expresses \(y|_u\) as a
product of at most \(\ell\) elements of \(C\cup C^{-1}\). Hence
\(\mathcal S_H(Y)\) is finite whenever \(C\) is finite.
\\
Assume now that every element of \(C\) has a section equal to \(1\).
The same holds for \(C^{-1}\): if \(c|_u=1\), then
\[
(c^{-1})|_{u\cdot c}=(c|_u)^{-1}=1.
\]
Let \(h=h_1\cdots h_\ell\), with \(h_j\in C\cup C^{-1}\). Choose
\(u_1\in A^*\) such that \(h_1|_{u_1}=1\). Then \(h|_{u_1}\) is a
product of at most \(\ell-1\) elements of \(C\cup C^{-1}\). Induction
on \(\ell\) gives \(u\in A^*\) such that \(h|_u=1\). Applying this to
every section of every \(y\in Y\) proves (ii). If \(H\) is self-similar, then \(A_i^*\) is \(H\)-invariant. The
restricted assertions follow by applying the preceding argument to
the induced self-similar action of \(H\) on \(A_i^*\).
\end{proof}

\begin{theorem}
\label{theo: finite-section-obstruction}
Let $H$ be self-similar and let $A=A_1\sqcup\cdots\sqcup A_k$ be its first-level orbit decomposition. Assume that, for every $i$, the graph
$\mathcal S_H^{(i)}(W)$ is finite and $1$ is reachable from every vertex. Let $K\leq H$ be a finitely generated self-similar subgroup
acting transitively on every $A_i$. Then $K$ is cyclic or non-free.
\end{theorem}
\begin{proof}
Suppose that \(K\) is a non-abelian free group and let
\(Y\subseteq K\subseteq H\) be a finite free basis of \(K\), with \(|Y|>1\). By Lemma~\ref{lem:finite-section-inheritance}, every $\mathcal S_H^{(i)}(Y)$ is finite and every vertex has a path to $1$. Lemma~\ref{lem: absorption-one} gives $\chi_H^{(i)}(Y)=1$ for every $i\in [1,k]$. Thus,  Proposition~\ref{pro:probabilistic-criterion}
implies that the natural map $\pi:F_Y\to K$ is not injective, a contradiction.
\end{proof}
The next result uses the full section graph.  The ambient group is
self-similar, but the subgroup is not assumed to be self-similar.
\begin{theorem}
\label{theo: level-transitive-obstruction}
Let $G=\langle X\rangle\leq\operatorname{Aut}(A^*)$ be
self-similar. Assume that
$\mathcal S_G(X)$ is finite and that $1$ is reachable from every
vertex. Every finitely generated level-transitive subgroup $H\leq G$
is cyclic or non-free.
\end{theorem}

\begin{proof}
Suppose that $H$ is a non-abelian free group with finite free basis
$W\subseteq H$ and $|W|>1$, and let
$\rho:F_W\to H$ be the natural isomorphism.
Lemma~\ref{lem:finite-section-inheritance} shows that
$\mathcal S_G(W)$ is finite and that $1$ is reachable from every
vertex. Lemma~\ref{lem: absorption-one} gives
$\chi_W(m)\to1$. Choose $m$ such that
$\chi_W(m)\geq1/|W|$. By
Lemma~\ref{lem:section-trivial-cycle}, there exist $u\in A^m$ and
$1\neq z\in F_W$ such that
\[
 u\cdot\rho(z)=u,\qquad \rho(z)|_u=1.
\]
Thus $z\in\ker(\sigma_u)$. Since $H$ is transitive on $A^m$,
Theorem~\ref{theo: local-global-criterion} shows that $W$ is not a
free basis of $H$, a contradiction.
\end{proof}

\begin{rem}
The proof uses transitivity only on the selected level $A^m$.  More
precisely, for a generating set $W$ as above it is enough that, for
some $m$, the subgroup $H$ be transitive on $A^m$ and
$\chi_{ W}(m)\geq1/|W|$.
\end{rem}
We now return to finite automata.  The transition graph of an automaton
projects onto the section graph of its states, so co-accessibility of an
identity state gives precisely the finite absorption hypothesis used
above.
\begin{cor}
\label{cor:coaccessible-identity}
Let \(\mathcal A\) be a finite invertible Mealy automaton over \(A\).
Assume that \(\mathcal A\) has a state \(e\) inducing the identity on
\(A^*\), reachable from every state. Let
\(H\leq G(\mathcal A)\) be finitely generated. If either
\begin{enumerate}[label=(\alph*)]
 \item \(H\) is self-similar and transitive on \(A\); or
 \item \(H\) is level-transitive on \(A^*\),
\end{enumerate}
then \(H\) is cyclic or non-free.
\end{cor}

\begin{proof}
Let \(X\) be the state set of \(\mathcal A\). If a path labelled by
\(u\in A^*\) starts at a state \(q\) and ends at a state \(p\), then
\(q|_u=p\). Hence the section graph
\(\mathcal S_{G(\mathcal A)}(X)\) is a quotient (some states of $X$ may be equal in $G(\mathcal{A})$) of the finite input
transition graph of \(\mathcal A\), obtained by forgetting the output
labels. Therefore \(\mathcal S_{G(\mathcal A)}(X)\) is finite.
Moreover, since \(e\) is reachable from every state and induces the
identity, \(1\) is reachable from every vertex of
\(\mathcal S_{G(\mathcal A)}(X)\).
\\
Now, if \(H\) is self-similar and transitive on \(A\), then
\(G(\mathcal A)\) also has a single first-level orbit, so
Theorem~\ref{theo: finite-section-obstruction} applies.
If \(H\) is level-transitive, apply
Theorem~\ref{theo: level-transitive-obstruction}.
\end{proof}
\begin{rem}
Theorem~\ref{theo: finite-section-obstruction} also applies when the
first-level action of \(H\) is not transitive. Indeed, if
\[
A=A_1\sqcup\cdots\sqcup A_k
\]
is the first-level orbit decomposition and \(H\) is self-similar, then
each \(A_i^*\) is \(H\)-invariant, and the action restricts to a
self-similar action on the rooted tree \(A_i^*\). Thus the argument may
be applied separately to every orbit, provided that \(1\) is reachable
in each restricted section graph \(\mathcal S_H^{(i)}(W)\). Global
co-accessibility in \(\mathcal S_H(W)\) does not imply this restricted
condition, since a path to \(1\) may use letters outside \(A_i\).
\end{rem}
\begin{figure}[ht]
\centering
\begin{tikzpicture}[
  >=Stealth,
  state/.style={circle,draw,minimum size=9mm,inner sep=1pt},
  every edge/.style={draw,->},
  every loop/.style={min distance=8mm}
]
\node[state] (a) at (0,1.65) {$a$};
\node[state] (b) at (0,-1.65) {$b$};
\node[state] (c) at (-2.55,0) {$c$};
\node[state] (e) at (3.15,0) {$e$};
\path
 (a) edge node[above left] {$1\mid0$} (c)
 (a) edge[bend right=20] node[left] {$0\mid1$} (b)
 (c) edge[loop left] node[left] {$0\mid1$} (c)
 (c) edge node[below left] {$1\mid0$} (b)
 (b) edge[bend right=20] node[right] {$\substack{0\mid0\\[-1pt]1\mid1}$} (a)
 (e) edge[loop right] node[right] {$\substack{0\mid0\\[-1pt]1\mid1}$} (e);
\end{tikzpicture}
\caption{The Aleshin automaton with an additional identity sink $e$.
The new state does not change the generated group, but it is not
reachable from the original component.}
\label{fig:aleshin-sink}
\end{figure}
\begin{rem}
Both the co-accessibility and transitivity assumptions are necessary.
Figure~\ref{fig:aleshin-sink} is obtained by adjoining an isolated
identity sink to the Aleshin automaton. The generated group remains
free of rank three \cite{Ale83,VV07}, although the identity state is no
longer reachable from the original component. Conversely, examples
in~\cite{fragile} show that a free automaton group may admit a
co-accessible identity state when the first-level action is not
transitive.
\end{rem}
The previous absorption argument also has a geometric consequence.
The relevant quantity is the number of vertices at which the generators
have non-trivial sections. If this number is negligible compared with
the size of the level, then the level Schreier graphs contain sets with
small boundary. More precisely, the same section-trivial edges used to produce local kernels can be grouped into components whose
boundaries are controlled by the complementary non-section-trivial
edges. For a finite graph $\Gamma=(V,E)$ and $S\subseteq V$, let
$\partial S$ be the set of edges with exactly one endpoint in $S$.
We write
\[
h(\Gamma):=
\min_{\substack{\varnothing\neq S\subseteq V\\ |S|\leq |V|/2}}
\frac{|\partial S|}{|S|}
\]
for the Cheeger constant of $\Gamma$. For the level Schreier graphs, we compute the Cheeger constant on the
underlying undirected multigraph: for each $u\in A^m$ and $x\in X$
there is one geometric edge joining $u$ to $u\cdot x$. This is the
undirected graph underlying the symmetric Schreier graph defined in
Section~2, since the $x^{-1}$-edge is the reverse of the corresponding
$x$-edge. Parallel edges and loops are kept.

\begin{prop}
\label{prop:activity-non-expansion}
Let $G=\langle X\rangle\leq\operatorname{Aut}(A^*)$, where $X$ is
finite, $d:=|A|\geq2$, and assume that $G$ acts level-transitively.
For $m\geq0$, put
\[
\alpha_X(m):=
\frac1{d^m}
\sum_{x\in X}
\bigl|\{u\in A^m:x|_u\neq1\}\bigr|.
\]
Let $\Gamma_m$ be the Schreier graph of the action of $G$ on $A^m$
with respect to $X$. Then
\[
h(\Gamma_{m+1})\leq2\alpha_X(m).
\]
In particular, if $\alpha_X(m)\to0$, then $h(\Gamma_m)\to0$.
Hence the level Schreier graphs do not form a family of expanders and
$G$ does not have Kazhdan's property~$(T)$.
\end{prop}

\begin{proof}
Let $E^{(m)}$ be the set of geometric $X$-edges of $\Gamma_m$
corresponding to pairs $(u,x)\in A^m\times X$ for which
$x|_u\neq1$. Thus \( |E^{(m)}|=d^m\alpha_X(m).\) Delete the edges in $E^{(m)}$ from $\Gamma_m$, and let
$\mathcal C_m$ be the set of connected components of the resulting multigraph. For $C\in\mathcal C_m$, let $b(C)$ be the
number of endpoints in $C$ of edges in $E^{(m)}$, counted with
multiplicity. Since the components partition $A^m$ and every edge has
two endpoints,
\[
\sum_{C\in\mathcal C_m}|C|=d^m,
\qquad
\sum_{C\in\mathcal C_m}b(C)=2|E^{(m)}|.
\]
Hence some $C\in\mathcal C_m$ satisfies
$b(C)/|C|\leq2\alpha_X(m)$. Fix $a\in A$ and put $S:=Ca\subseteq A^{m+1}$. We claim that
$|\partial S|\leq b(C)$. Indeed, let
$ua\xrightarrow{x}vb$ be an edge of $\Gamma_{m+1}$ with exactly one
endpoint in $S$. Its projection to level $m$ is the edge
$u\xrightarrow{x}v$ of $\Gamma_m$. If this projected edge were
section-trivial, then $x|_u=1$, so $b=a$ and
\[
(ua)\cdot x=(u\cdot x)a=va.
\]
Moreover, $u$ and $v$ belong to the same component of the
section-trivial subgraph. Hence $ua\in S$ if and only if $va\in S$,
contrary to the choice of the edge. Thus every edge of $\partial S$
projects to an edge of $E^{(m)}$ having an endpoint in $C$. For each endpoint in $C$ of an edge of $E^{(m)}$, there is at most
one edge of $\partial S$ projecting to it: this is immediate at the
initial endpoint, while at the terminal endpoint it follows from the
fact that the section acts bijectively on $A$. Therefore
$|\partial S|\leq b(C)$. Finally, since $d\geq2$,
$|S|=|C|\leq d^m\leq |A^{m+1}|/2$. Hence
\[
h(\Gamma_{m+1})
\leq\frac{|\partial S|}{|S|}
\leq\frac{b(C)}{|C|}
\leq2\alpha_X(m).
\]
Thus $\alpha_X(m)\to0$ implies $h(\Gamma_m)\to0$. If $G$ had property~\emph{(T)}, the finite transitive $G$-sets $A^m$
would have a uniform spectral gap, and hence a uniform positive
Cheeger constant, a contradiction.
\end{proof}

\begin{rem}
Spectral and growth properties of Schreier graphs of automaton groups
have been studied, for instance, in
\cite{GrigorchukSunicSchreier,BondarenkoSchreier}; the former
investigates, in particular, the decay of the spectral gap, while the
latter treats Schreier graphs associated with automata of polynomial
activity. Proposition~\ref{prop:activity-non-expansion} gives a
complementary criterion for the level Schreier graphs: it only
requires that the proportion of vertices carrying a non-trivial
section tend to zero. In particular, it applies to finite generating
sets of polynomial activity and, independently, to finite automata
with a co-accessible identity state.
\end{rem}

The preceding criterion is an immediate consequence of absorption. Indeed, with the above notation we have \(\alpha_X(m)=|X|\bigl(1-\chi_X(m)\bigr).\)

\begin{cor}\label{cor:coaccessible-no-T}
Let $\mathcal A$ be a finite invertible Mealy automaton over an alphabet
$A$ with $|A|\geq2$. Assume that $\mathcal A$ has a co-accessible
identity state and that $G(\mathcal A)$ acts level-transitively. Then
$h(\Gamma_m)\to0$ for the level Schreier graphs of $G(\mathcal A)$.
In particular, they do not form a family of expanders and
$G(\mathcal A)$ does not have Kazhdan's property~\emph{(T)}.
\end{cor}

\begin{proof}
By Lemma~\ref{lem: absorption-one}, $\chi_X(m)\to1$, hence
$\alpha_X(m)=|X|(1-\chi_X(m))\to0$. The conclusion follows from
Proposition~\ref{prop:activity-non-expansion}.
\end{proof}

\begin{cor}\label{cor:polynomial-activity-no-T}
Let $G=\langle X\rangle\leq\Aut(A^\ast)$ be a finitely generated
self-similar group over an alphabet $A$ with $|A|\geq2$. Assume that
$G$ acts level-transitively and that every $x\in X$ has polynomial
activity. Then $h(\Gamma_m)\to0$ for the level Schreier graphs of $G$.
In particular, they do not form a family of expanders and $G$ does not
have Kazhdan's property~\emph{(T)}.
\end{cor}

\begin{proof}
Polynomial activity gives
$|\{u\in A^m:x|_u\neq1\}|=O(m^{k_x})$ for every $x\in X$. Since $X$
is finite, $\alpha_X(m)=O(m^k/|A|^m)\to0$ for some $k\geq0$, and the
conclusion follows from Proposition~\ref{prop:activity-non-expansion}.
\end{proof}

\section{Effective and quantitative forms}
\label{sec: effective forms}

This section gives effective versions of the preceding obstruction.
For finite section graphs, the absorption probabilities can be computed
exactly.  Once the inequality in Proposition~\ref{pro:probabilistic-criterion}
holds at some finite depth, the proof constructs an explicit relation;
keeping track of the lengths gives a coarse upper bound.

\subsection{Absorption probabilities}

Let $H=\langle W\rangle\leq\Aut(A^*)$ be finitely generated and
self-similar, let $\rho:F_W\to H$ be the natural epimorphism, and write
$A=A_1\sqcup\cdots\sqcup A_k$ for the decomposition into $H$-orbits.
Assume that every restricted section graph
$\mathcal S_H^{(i)}(W)$ is finite.  For a vertex $v$ of this graph, let
$h_i(v)$ be the probability that the random walk whose labels are
chosen uniformly in $A_i$ ever reaches the identity. The following facts are standard facts in the theory of Markov chains. 

\begin{prop}
\label{prop: absorption-criterion}
For every \(i\), the absorption probabilities \(h_i(v)\) are rational
and can be computed from the finite section graph. Moreover,
\[
\chi_H^{(i)}(W)
=
\frac1{|W|}
\sum_{w\in W}h_i(w).
\]
If \(\sum_{w\in W}h_i(w)>1\) for \(i\in [1,k]\) then the natural epimorphism \(\rho:F_W\to H\) is not injective.
\\
If the self-similar action and the finite section graphs are given
effectively, then a non-trivial element of \(\ker(\rho)\) can be
constructed.
\end{prop}
\begin{proof}
Consider the vertices that are co-accessible from the identity:
\[
V_i:=\{v\neq1:1\text{ is reachable from }v\},
\qquad
h_i:=(h_i(v))_{v\in V_i}.
\]
For \(v\in V_i\),
\[
h_i(v)=\frac1{|A_i|}\sum_{a\in A_i}h_i(v|_a),
\]
where \(h_i(1)=1\) and \(h_i(v|_a)=0\) whenever \(1\) is not reachable
from \(v|_a\). Hence, \(h_i=Q_i h_i+b_i,\)
where \(Q_i\) records the transitions inside \(V_i\) and \(b_i\) the
transitions from \(V_i\) to \(1\). Thus
\[
(I-Q_i)h_i=b_i.
\]
Since every vertex of \(V_i\) reaches \(1\), one has \(Q_i^n\to0\).
Therefore \(I-Q_i\) is invertible, and
\[
h_i=(I-Q_i)^{-1}b_i.
\]
The entries of \(Q_i\) and \(b_i\) are rational, so every \(h_i(v)\) is
rational and computable. Moreover, \(\chi_H^{(i)}(W)
=
\frac1{|W|}\sum_{w\in W}h_i(w).\)
Now, if \(\sum_{w\in W}h_i(w)>1\) for each $i\in [1,k]$, then, there exists \(m_i\) such that
\[
\sum_{w\in W}
|\{u\in A_i^{m_i}:w|_u=1\}|
>
|A_i|^{m_i}.
\]
Lemma~\ref{lem:section-trivial-cycle} gives
\(u_i\in A_i^{m_i}\) and \(1\neq z_i\in F_W\) such that \(u_i\cdot\rho(z_i)=u_i, \,\rho(z_i)|_{u_i}=1.\)
Thus \(1\neq z_i\in\ker(\sigma_{u_i})\) for every first-level orbit.
Corollary~\ref{cor: first level reduction} yields \(\ker(\rho)\neq1\).
\\
If the action and the section graphs are given effectively, the
integers \(m_i\), the cycles of
Lemma~\ref{lem:section-trivial-cycle}, and the corresponding local kernel elements \(z_i\) can be computed. Applying Lemma~\ref{lem: local-kernel-patching} and Proposition~\ref{prop: level-kernel} then produces a non-trivial element of \(\ker(\rho)\).
\end{proof}

\subsection{A bound for the shortest relation}

We now keep track of lengths in the local-to-global construction. Let
\(G\leq\Aut(A^*)\) be self-similar and let \(H=\langle W\rangle\leq G\),
where \(W\) is finite and \(|W|\geq2\). Let
\(\rho:F_W\twoheadrightarrow H\) be the natural epimorphism, and denote by
\(|\cdot|_W\) the word length in \(F_W\). Put
\[
  \operatorname{rel}_W(H)
  :=\min\{|z|_W:1\neq z\in\ker(\rho)\},
\]
with value \(\infty\) if \(\rho\) is injective. We denote by \(\mathcal L(q)\) \emph{Landau's function}, the maximal order of
an element of \(\operatorname{Sym}(q)\). Recall that
\begin{equation}\label{eq:landau}
  \mathcal L(q)=\exp\bigl((1+o(1))\sqrt{q\log q}\,\bigr)
\end{equation}
as \(q\to\infty\) \cite{Landau}. For \(g\in F_W\) and \(n\geq1\), let
\(o_n(g)\) be the order of the permutation induced by \(\rho(g)\) on \(A^n\).
Then \(g^{o_n(g)}\in P_n\) and \(o_n(g)\leq\mathcal L(|A|^n)\).
 For \(S\subseteq A^n\), an \emph{\(S\)-kernel} is an element
\(1\neq w\in P_n\cap\bigcap_{v\in S}\ker(\sigma_v)\). By
Proposition~\ref{prop: level-kernel}, an \(A^n\)-kernel is precisely a
non-trivial element of \(\ker(\rho)\). We first merge two \(S\)-kernels while
controlling their length.

\begin{lemma}\label{lem:quantitative-merge}
Assume \(\operatorname{rk}(F_W)\geq2\), fix \(n\geq1\), and put
\(q:=|A|^n\). Let \(S_1,S_2\subseteq A^n\) and let \(w_i\) be an
\(S_i\)-kernel for \(i=1,2\). Then there exists an
\((S_1\cup S_2)\)-kernel \(w\) such that
\[
  |w|_W\leq2|w_1|_W+2|w_2|_W+4\mathcal L(q).
\]
Moreover, \(w\) can be computed from \(w_1\) and \(w_2\).
\end{lemma}

\begin{proof}
By Lemma~\ref{lem: non trivial commutators}, choose
\(c\in\{1\}\cup W\) such that
\([c^m w_1c^{-m},w_2]\neq1\) for every \(m\in\mathbb Z\setminus\{0\}\).
Put \(E:=o_n(c)\) and \(w:=[c^Ew_1c^{-E},w_2]\neq1\). Since \(c^E,w_1,w_2\in P_n\), we have \(w\in P_n\). If \(v\in S_1\), then
\(\sigma_v(c^Ew_1c^{-E})=1\), while if \(v\in S_2\), then
\(\sigma_v(w_2)=1\). Hence \(\sigma_v(w)=1\) for every
\(v\in S_1\cup S_2\), so \(w\) is an \((S_1\cup S_2)\)-kernel. Finally, \(|c^Ew_1c^{-E}|_W\leq|w_1|_W+2E\), with
\(E\leq\mathcal L(q)\), and the stated bound follows. The choice of \(c\) is
effective by Lemma~\ref{lem: non trivial commutators}, while \(E\) is computed
from the finite action on \(A^n\).
\end{proof}

\begin{prop}\label{prop:quantitative-patching}
Assume that \(H\) acts transitively on \(A^n\) and put \(q:=|A|^n\). If, for
some \(u\in A^n\), there exists \(1\neq z\in\ker(\sigma_u)\) with
\(|z|_W\leq\ell\), then
\[
  \operatorname{rel}_W(H)\leq4q^2\mathcal L(q)(\ell+2q).
\]
\end{prop}

\begin{proof}
For each \(v\in A^n\), choose \(t_v\in F_W\) with
\(v\cdot\rho(t_v)=u\) and \(|t_v|_W\leq q-1\); such a word is obtained from a
spanning tree of the symmetric Schreier graph on \(A^n\). By the transport
computation in Lemma~\ref{lem: local-kernel-patching},
\(z_v:=t_vzt_v^{-1}\) is non-trivial, belongs to \(\ker(\sigma_v)\), and
satisfies \(|z_v|_W\leq\ell+2q-2\). Set \(y_v:=z_v^{\,o_n(z_v)}\). Since \(F_W\) is torsion-free,
\(y_v\neq1\); moreover \(y_v\in P_n\) and \(\sigma_v(y_v)=1\). Thus \(y_v\)
is a \(\{v\}\)-kernel and
\(|y_v|_W\leq\mathcal L(q)(\ell+2q-2)=:M_0\). We merge these \(q\) singleton kernels recursively along a balanced binary
tree. Put \(C:=\frac43\mathcal L(q)\). We claim that, for every non-empty
\(S\subseteq A^n\), the \(S\)-kernel \(w_S\) obtained in this way satisfies
\[
  |w_S|_W
  \leq4^{\lceil\log_2|S|\rceil}(M_0+C)-C.
\]
For \(|S|=1\) this follows from the definition of \(M_0\). Let
\(|S|=k\geq2\), and split \(S=S_1\sqcup S_2\) with
\(|S_1|=\lceil k/2\rceil\) and \(|S_2|=\lfloor k/2\rfloor\). Then
\(\lceil\log_2|S_i|\rceil\leq\lceil\log_2k\rceil-1\). By
Lemma~\ref{lem:quantitative-merge} and the induction hypothesis,
\[
\begin{aligned}
 |w_S|_W
 &\leq2|w_{S_1}|_W+2|w_{S_2}|_W+4\mathcal L(q)\\
 &\leq4^{\lceil\log_2k\rceil}(M_0+C)-4C+4\mathcal L(q)\\
 &=4^{\lceil\log_2k\rceil}(M_0+C)-C,
\end{aligned}
\]
since \(4\mathcal L(q)=3C\). Taking \(S=A^n\) and using
\(4^{\lceil\log_2q\rceil}\leq4q^2\), we obtain an \(A^n\)-kernel \(w\) with
\[
  |w|_W
  \leq4q^2(M_0+C)
  \leq4q^2\mathcal L(q)(\ell+2q).
\]
By Proposition~\ref{prop: level-kernel}, \(1\neq w\in\ker(\rho)\).
\end{proof}
Assume now that \(\mathcal S_G(W)\) is finite and that \(1\) is reachable from
every vertex. Let \(d_{\mathcal S}(v,v')\) be the length of a shortest
directed walk from \(v\) to \(v'\), with value \(\infty\) if no such walk
exists, and put
\(L:=\max_{v\in V(\mathcal S_G(W))}d_{\mathcal S}(v,1)\),
\(d:=|A|\), and \(r:=|W|\).
If \(L=0\), then every element of \(W\) is trivial and
\(\operatorname{rel}_W(H)=1\). Assume \(L\geq1\), and set
\[
  s_0:=\min\left\{s\geq1:(1-d^{-L})^s<\frac{r-1}{r}\right\},
  \qquad m_0:=Ls_0.
\]

\begin{lemma}\label{lem:quantitative-absorption}
For every vertex \(v\) of \(\mathcal S_G(W)\) and every \(m\geq1\), the
probability that a uniformly chosen word in \(A^m\) sends \(v\) to \(1\) is at
least \(1-(1-d^{-L})^{\lfloor m/L\rfloor}\). Consequently,
\(|T(m_0)|>d^{m_0}\).
\end{lemma}

\begin{proof}
From every vertex there is a word of length at most \(L\) leading to \(1\).
Hence, conditionally on not having reached \(1\), the probability of reaching
it during the next \(L\) steps is at least \(d^{-L}\). Thus the probability of
avoiding \(1\) for \(sL\) steps is at most \((1-d^{-L})^s\). For every \(w\in W\), applying the preceding estimate with
\(m=m_0=Ls_0\) gives
\[
  \frac{|\{u\in A^{m_0}:\rho(w)|_u=1\}|}{d^{m_0}}
  \geq 1-(1-d^{-L})^{s_0}
  >\frac1r.
\]
Summing over \(w\in W\), and recalling that
\(T(m_0)=\bigsqcup_{w\in W}T_w(m_0)\), we obtain \(  |T(m_0)|
  >d^{m_0}.\)
\end{proof}

Together with Lemma~\ref{lem:section-trivial-cycle}, the preceding estimate
gives a local kernel of controlled length. Proposition~\ref{prop:quantitative-patching}
then gives a global relation.

\begin{theorem}\label{thm:quantitative}
Let \(G\leq\operatorname{Aut}(A^*)\) be self-similar and let
\(H=\langle W\rangle\leq G\), where \(W\) is finite and \(|W|\geq2\).
Assume that \(\mathcal S_G(W)\) is finite and that \(1\) is reachable from
every vertex. Let \(d,r,L\) be as above. If \(L=0\), then
\(\operatorname{rel}_W(H)=1\). Suppose \(L\geq1\), and let \(s_0,m_0\) be as
above. Then:
\begin{enumerate}
\item[(i)] If \(H\) is level-transitive and \(q=d^{m_0}\), then
\[
  \operatorname{rel}_W(H)\leq12q^3\mathcal L(q).
\]

\item[(ii)] Suppose that \(H\) is self-similar and transitive on \(A\).
For \(a\in A\) and \(x\in W\), choose
\(\lambda(a,x)\in F_W\) representing \(\rho(x)|_a\), and set
\(\lambda(a,x^{-1}):=\lambda(a\cdot\rho(x)^{-1},x)^{-1}\). Let
\[
  C:=\max\Bigl(
  \{1\}\cup
  \{|\lambda(a,x^\varepsilon)|_W:
  a\in A,\ x\in W,\ \varepsilon\in\{\pm1\}\}
  \Bigr).
\]
Then
\[
  \operatorname{rel}_W(H)
  \leq4d^2\mathcal L(d)
  \bigl(C^{m_0-1}d^{m_0}+2d\bigr).
\]
If every first-level section of an element of
\(W\cup W^{-1}\) belongs to \(W\cup W^{-1}\cup\{1\}\), then one may take
\(C=1\).
\end{enumerate}
\end{theorem}

\begin{proof}
If \(L=0\), then every vertex of \(\mathcal S_G(W)\) equals \(1\), so every
element of \(W\) is trivial. Assume \(L\geq1\). By
Lemma~\ref{lem:quantitative-absorption}, \(|T(m_0)|>d^{m_0}\). Hence the
undirected multigraph \(\Gamma_{m_0}\) from the proof of
Lemma~\ref{lem:section-trivial-cycle} has more edges than vertices. Some
connected component therefore contains a simple cycle, where a loop or a pair
of parallel edges is allowed. Its length is at most \(d^{m_0}\). Repeating the
labelling argument of Lemma~\ref{lem:section-trivial-cycle} gives
\(u\in A^{m_0}\) and \(z\in F_W\) such that
\begin{equation}\label{eq:localkernel}
  1\neq z,\qquad
  u\cdot\rho(z)=u,\qquad
  \rho(z)|_u=1,\qquad
  |z|_W\leq d^{m_0}.
\end{equation}
Thus \(1\neq z\in\ker(\sigma_u)\).

For \textup{(i)}, apply Proposition~\ref{prop:quantitative-patching} with
\(n=m_0\), \(q=d^{m_0}\), and \(\ell=q\). Then
\(\operatorname{rel}_W(H)\leq4q^2\mathcal L(q)(q+2q)
=12q^3\mathcal L(q)\).

For \textup{(ii)}, write \(u=a_1\cdots a_{m_0}\) and put
\(g_0:=\rho(z)\) and
\(g_i:=\rho(z)|_{a_1\cdots a_i}\) for \(1\leq i\leq m_0\). Then
\(g_{m_0}=1\). If \(g_0=1\), then
\(z\in\ker(\rho)\setminus\{1\}\) and
\(\operatorname{rel}_W(H)\leq d^{m_0}\), which is bounded by the claimed
right-hand side. Assume \(g_0\neq1\), and let \(j\geq1\) be minimal with \(g_j=1\).
Put \(z_0:=z\) and, for \(1\leq i<j\), define
\(z_i:=\widehat\sigma_{a_i}(z_{i-1})\). By
Lemma~\ref{lem: lifted section map}, \(\rho(z_i)=g_i\). Moreover
\(|z_i|_W\leq C|z_{i-1}|_W\), hence
\(|z_i|_W\leq C^i|z|_W\). Set \(z':=z_{j-1}\). By minimality of \(j\),
\(\rho(z')=g_{j-1}\neq1\), hence \(z'\neq1\). Since \(\rho(z)\) fixes \(u\),
the section \(g_{j-1}\) fixes \(a_j\), and
\(\sigma_{a_j}(z')=g_{j-1}|_{a_j}=g_j=1\). Thus
\(1\neq z'\in\ker(\sigma_{a_j})\) and
\(|z'|_W\leq C^{m_0-1}d^{m_0}\). Applying Proposition~\ref{prop:quantitative-patching} to the transitive
first-level action, with \(n=1\), \(q=d\), and
\(\ell=C^{m_0-1}d^{m_0}\), gives the claimed bound. Finally, if every
first-level section of an element of \(W\cup W^{-1}\) lies in
\(W\cup W^{-1}\cup\{1\}\), the representatives
\(\lambda(a,x^\varepsilon)\) may be chosen of length at most \(1\), so \(C=1\).
\end{proof}

\begin{rem}
The bounds above are not optimal. We record their asymptotic behaviour as
\(L\to\infty\), with \(d=|A|\geq2\) and \(r=|W|\) fixed. Put
\(\alpha:=\log\frac{r}{r-1}\). From the definition of \(s_0\),
\(s_0=\alpha d^L+O(1)\), and hence
\(m_0=\alpha Ld^L+O(L)\) and
\(q_0:=d^{m_0}
=\exp(\alpha\log(d)\,Ld^L+O(L))\).
By \eqref{eq:landau}, Theorem~\ref{thm:quantitative}\textup{(i)} gives
\[
  \log\operatorname{rel}_W(H)
  \leq(1+o(1))\sqrt{q_0\log q_0},
\]
and therefore
\(\log\log\operatorname{rel}_W(H)
\leq\frac12\log q_0+O(\log\log q_0)\).
For fixed \(C\), Theorem~\ref{thm:quantitative}\textup{(ii)} gives
\(\operatorname{rel}_W(H)=\exp(O(Ld^L))\); if \(C=1\), the bound is
linear in \(q_0\), up to a constant depending only on \(d\).
\end{rem}

\section{Section collisions and bireversibility}
\label{sec:reversible-not-bireversible}

The second obstruction comes from collisions of sections. We first show
that a collision of a surjective section map produces a local kernel;
later we obtain the same conclusion from a Stallings fold when the
states form a free basis. Let \(H\leq\operatorname{Aut}(A^*)\) be
self-similar and, for \(u\in A^*\), write
\[
\operatorname{sec}_u(g):=g|_u.
\]
Thus \(\operatorname{sec}_u:H\to H\) is a map of sets. Its restriction
to \(\operatorname{Stab}_H(u)\) is the canonical section
homomorphism. Let \(B\subseteq A\) be an \(H\)-orbit. The action of \(H\) is
\emph{section-surjective on the subtree \(B^*\)} if
\(\operatorname{sec}_u\) is surjective for every \(u\in B^*\). It is
\emph{section-bijective on \(B^*\)} if every \(\operatorname{sec}_u\),
\(u\in B^*\), is bijective. We write
\[
\theta_u
:=
\operatorname{sec}_u|_{\operatorname{Stab}_H(u)}
:
\operatorname{Stab}_H(u)\longrightarrow H
\]
for the section homomorphism on the vertex stabilizer. The terminology is intrinsic to the self-similar action. By
Proposition~\ref{prop:reversible-surjective}, every group generated by
a finite reversible invertible Mealy automaton acts
section-surjectively on each first-level orbit.

\subsection{From a collision of sections to a local kernel}

A collision for \(g\mapsto g|_t\) does not necessarily occur inside
\(\operatorname{Stab}_H(t)\), so it does not directly produce an
element of \(\ker(\theta_t)\). The next lemma shows that, after moving
along the finite orbit of \(t\), one obtains a genuine local kernel.

\begin{lemma}
\label{lem:section-collision}
Let \(H\) be a finitely generated non-abelian free group acting
self-similarly on \(A^*\), and let \(t\in A^m\). If
\(\operatorname{sec}_t:H\to H\) is surjective but not injective, then
there exists \(v\in t\cdot H\) such that
\[
\ker(\theta_v)\neq1.
\]
\end{lemma}

\begin{proof}
Consider the orbit \(O:=t\cdot H\) and put \(N:=|O|\) Suppose, for a contradiction, that
\(\theta_v\) is injective for every \(v\in O\). For each \(v\in O\), choose \(r_v\in H\) such that
\(t\cdot r_v=v\), and put
\[
c_v:=r_v|_t,
\qquad
L_v:=\theta_v(\operatorname{Stab}_H(v)).
\]
The fiber of the action map over \(v\) is precisely the left coset
\(r_v\operatorname{Stab}_H(v)\). Indeed,
\[
t\cdot g=v
\iff
r_v^{-1}g\in\operatorname{Stab}_H(v).
\] Hence \(g|_t
=
r_v|_t\,k|_{t\cdot r_v}
=
c_v\theta_v(k).\)
Thus the image of the fiber \(\{g\in H:t\cdot g=v\}\) under
\(\operatorname{sec}_t\) is \(c_vL_v\). Since
\(\operatorname{sec}_t\) is surjective,
\begin{equation}
\label{eq:section-coset-cover}
H=\bigcup_{v\in O}c_vL_v.
\end{equation}
For \(k\in\operatorname{Stab}_H(v)\), the element
\(r_vkr_v^{-1}\) belongs to \(\operatorname{Stab}_H(t)\), and \(\theta_t(r_vkr_v^{-1})
=
c_v\theta_v(k)c_v^{-1}.\) Conjugation by \(r_v\) maps
\(\operatorname{Stab}_H(v)\) onto
\(\operatorname{Stab}_H(t)\); therefore \(L_t=c_vL_vc_v^{-1}.\) Hence the subgroups \(L_v\), \(v\in O\), are conjugate.
\\
By a theorem of B.~H.~Neumann~\cite{NeumannCover}, whenever a group is
covered by finitely many cosets of subgroups, at least one of the
subgroups has finite index. Hence by \eqref{eq:section-coset-cover}, at least one \(L_v\) has finite index in \(H\). Consequently,
all \(L_v\) have the same finite index; write \([H:L_v]=d.\) Let \(q:=\operatorname{rk}(H)\geq2\). Since the action of \(H\) on
\(O\) is transitive,
\[
[H:\operatorname{Stab}_H(v)]=N.
\]
By Schreier's formula, \(\operatorname{rk}(\operatorname{Stab}_H(v))
=
1+N(q-1).\)
The injectivity of \(\theta_v\) gives
\(L_v\simeq\operatorname{Stab}_H(v)\), whereas
\([H:L_v]=d\) gives \(\operatorname{rk}(L_v)=1+d(q-1).\) Since \(q\geq2\), it follows that \(d=N\). Set
\[
K:=\bigcap_{v\in O}\operatorname{Core}_H(L_v).
\]
The subgroup \(K\) is normal and of finite index in \(H\), and
\(K\leq L_v\) for every \(v\in O\). Let \(Q:=H/K\). Each coset
\(c_vL_v\) projects to a subset of \(Q\) of cardinality
\[
\frac{|Q|}{[H:L_v]}
=
\frac{|Q|}{N}.
\]
These \(N\) subsets cover \(Q\) by
\eqref{eq:section-coset-cover}, and the sum of their cardinalities is
\(|Q|\). Hence they are pairwise disjoint. The original cosets
\(c_vL_v\) are therefore pairwise disjoint as well. The fibers of \(g\mapsto t\cdot g\) partition \(H\). On the fiber
\(r_v\operatorname{Stab}_H(v)\), the map
\(\operatorname{sec}_t\) is injective because \(\theta_v\) is
injective, and its image is \(c_vL_v\). Since these images are pairwise
disjoint, \(\operatorname{sec}_t\) is injective on \(H\), a
contradiction.
\end{proof}
\begin{theorem}
\label{theo:section-surjective-collision}
Let \(H\leq\operatorname{Aut}(A^*)\) be a finitely generated
self-similar group, and let
\(A=A_1\sqcup\cdots\sqcup A_k\) be the decomposition into \(H\)-orbits
on the first level. Assume that, for every \(j\in\{1,\ldots,k\}\),
there exists \(t_j\in A_j^*\) such that
\[
\operatorname{sec}_{t_j}:H\to H
\]
is surjective but not injective. Then \(H\) is cyclic or not free.
\end{theorem}

\begin{proof}
Suppose, for a contradiction, that \(H\) is a non-abelian free group.
For each \(j\), Lemma~\ref{lem:section-collision}, applied to \(t_j\),
gives a vertex
\(v_j\in t_j\cdot H\subseteq A_j^{|t_j|}\)
such that \(\ker(\theta_{v_j})\neq1\).

Let \(W\) be a free basis of \(H\), and let
\(\rho:F_W\to H\) be the natural isomorphism. For each \(j\), choose
\(1\neq g_j\in\ker(\theta_{v_j})\) and put
\(z_j:=\rho^{-1}(g_j)\). Then \(z_j\neq1\),
\(\rho(z_j)\in\operatorname{Stab}_H(v_j)\), and
\(\rho(z_j)|_{v_j}=1\). Hence
\(z_j\in\ker(\sigma_{v_j})\).

Thus, for every first-level orbit \(A_j\), there exists
\(v_j\in A_j^{|v_j|}\) such that
\(\ker(\sigma_{v_j})\neq1\). By
Corollary~\ref{cor: first level reduction}, \(W\) is not a free basis
of \(H\), a contradiction.
\end{proof}

\begin{cor}
\label{cor:intrinsic-bijective-orbit}
Let \(H\) be a finitely generated self-similar non-abelian free group.
If the action is section-surjective on every first-level orbit, then
there exists a first-level orbit \(B\subseteq A\) on which the action
is section-bijective.
\end{cor}

\begin{proof}
Suppose that the action is not section-bijective on any first-level
orbit. Since it is section-surjective on every such orbit, for each
first-level orbit \(A_j\) there exists \(t_j\in A_j^*\) such that
\(\operatorname{sec}_{t_j}\) is not injective.
Theorem~\ref{theo:section-surjective-collision} then implies that \(H\)
is cyclic or not free, a contradiction.
\end{proof}

\subsection{Application to reversible automata}

We now translate the intrinsic criterion into automata-theoretic
language. Reversibility yields surjective section maps, whereas
non-bireversibility produces collisions. Thus the hypotheses of
Theorem~\ref{theo:section-surjective-collision} can be verified
directly on the automaton. Let \(\mathcal A=(X,A,\delta,\lambda)\) be a finite invertible Mealy
automaton, with \(x|_a:=\delta(x,a)\) and
\(a\cdot x:=\lambda(x,a)\). For \(m\geq0\), we identify the \(m\)-th power of the dual automaton
with the Mealy automaton whose state set is \(A^m\), alphabet \(X\),
and transitions
\[
v\xrightarrow{\,q\mid q|_v\,}v\cdot q.
\]

\begin{prop}
\label{prop:reversible-surjective}
If \(\mathcal A\) is reversible, then \(\operatorname{sec}_u:G(\mathcal A)\longrightarrow G(\mathcal A)\) is surjective for every \(u\in A^*\).
\end{prop}
\begin{proof}
Fix \(m\geq0\) and \(x\in X\). By reversibility, for every
\(v\in A^m\) there exists a unique state \(q_v\in X\) such that
\(q_v|_v=x\). Thus, in the \(m\)-th power of the dual, each vertex
\(v\) has a unique outgoing transition with output \(x\), namely
\(v\xrightarrow{\,q_v\mid x\,}v\cdot q_v.\)
Now fix \(u\in A^m\). To construct an element having section
\(x^{-1}\) at \(u\), set \(v_0:=u\) and define recursively
\[
v_{i+1}:=v_i\cdot q_{v_i}.
\]
Thus \(v_i\xrightarrow{\,q_{v_i}\mid x\,}v_{i+1}\) is precisely the unique outgoing edge labelled by \(x\) in the output from \(v_i\). Since \(A^m\) is finite, there exist \(k\geq0\) and \(p\geq1\) such
that \(v_k=v_{k+p};\) see Figure~\ref{fig:reversible-surjective}.
\begin{figure}[t]
\centering
\begin{tikzpicture}[
    >=Stealth,
    vertex/.style={circle,draw,minimum size=8mm,inner sep=1pt},
    forward/.style={->,thick},
    backward/.style={->,thick,dashed}
]
\node[vertex] (v0) at (0,0) {\(v_0=u\)};
\node[vertex] (v1) at (2.2,0) {\(v_1\)};
\node (dots1) at (4,0) {\(\cdots\)};
\node[vertex] (vk) at (5.5,0) {\(v_k\)};

\node[vertex] (vk1) at (7.5,2) {\(v_{k+1}\)};
\node (dots2) at (9.2,2) {\(\cdots\)};
\node[vertex] (vkp) at (11,0) {\(v_{k+p-1}\)};

\draw[forward]
(v0) -- node[above] {\(q_{v_0}\mid x\)} (v1);

\draw[forward] (v1) -- (dots1);
\draw[forward] (dots1) -- (vk);

\draw[forward]
(vk) to[bend left=12]
node[above left] {\(q_{v_k}\mid x\)}
(vk1);

\draw[forward] (vk1) -- (dots2);
\draw[forward] (dots2) -- (vkp);

\draw[forward]
(vkp) to[bend left=25]
node[below] {\(q_{v_{k+p-1}}\mid x\)}
(vk);

\draw[backward]
(vk1) to[bend left=10]
node[right] {\(q_{v_k}^{-1}\mid x^{-1}\)}
(vk);

\draw[backward]
(vk) to[bend left=10]
node[below] {\(q_{v_{k+p-1}}^{-1}\mid x^{-1}\)}
(vkp);

\node[align=center] at (8.4,-2)
{\(\text{the cycle is traversed backwards through }k+1\text{ edges}\)};
\end{tikzpicture}
\caption{The path from \(u\) to \(v_k\) contributes the section
\(x^k\), whereas \(k+1\) backwards edges in the cycle contribute
\(x^{-(k+1)}\).}
\label{fig:reversible-surjective}
\end{figure}
Every edge of the resulting cycle can be traversed backwards. By the
inverse-section identity,
\[
q_{v_i}^{-1}|_{v_{i+1}}
=
(q_{v_i}|_{v_i})^{-1}
=
x^{-1}.
\]
Starting at \(v_k\), traverse the cycle backwards through exactly
\(k+1\) edges, repeating the cycle when necessary, and let \(z\) be
the product of their inverse labels in traversal order. Iterating the
product identity gives \(z|_{v_k}=x^{-(k+1)}.\)
Put \(h_0:=q_{v_0}\cdots q_{v_{k-1}},\) with \(h_0=1\) when \(k=0\). Then \(u\cdot h_0=v_k,\) \(h_0|_u=x^k.\) Consequently,
\[
(h_0z)|_u
=
h_0|_u\,z|_{u\cdot h_0}
=
x^kx^{-(k+1)}
=
x^{-1}.
\]
Thus every element of \(X\cup X^{-1}\) occurs as a section at every
vertex of \(A^*\). Let \(g=s_1\cdots s_\ell\in G(\mathcal A)\), where
\(s_i\in X\cup X^{-1}\). Starting with \(u_0:=u\), choose
\(h_i\in G(\mathcal A)\) such that
\[
h_i|_{u_{i-1}}=s_i,
\qquad
u_i:=u_{i-1}\cdot h_i.
\]
Iterating the product identity yields \((h_1\cdots h_\ell)|_u
=
s_1\cdots s_\ell
=
g.\) Hence \(\operatorname{sec}_u\) is surjective.
\end{proof}

\begin{cor}
\label{cor:reversible-dual}
Let \(\mathcal A\) be a finite reduced reversible invertible Mealy
automaton. If every connected component of the dual automaton
\(\partial\mathcal A\) is non-bireversible, then
\(G(\mathcal A)\) is cyclic or not free.
\end{cor}

\begin{proof}
The connected components of \(\partial\mathcal A\) are the orbits
\(A_1,\ldots,A_k\) of \(G(\mathcal A)\) on \(A\). By
Proposition~\ref{prop:reversible-surjective}, the action is
section-surjective on every orbit subtree \(A_j^*\). Fix \(j\), and regard \(A_j\) as the induced Mealy automaton described
in Section~\ref{sec: preliminaries}. For the induced component on \(A_j\), the transition maps are
\(a\mapsto a\cdot x\) (\(x\in X\)), while the output maps are
\(x\mapsto x|_a\) (\(a\in A_j\)). The former are permutations because
\(\mathcal A\) is invertible, and the latter are permutations because
\(\mathcal A\) is reversible. Thus the component is reversible and
invertible, and since it is not bireversible, its inverse is not reversible. Hence there exist distinct states
\(a_1^{-1}\neq a_2^{-1}\), an input \(y\in X\), and transitions
\[
a_i^{-1}\xrightarrow{\,y\mid x_i\,}b^{-1}
\qquad(i=1,2)
\]
in the inverse component. Equivalently, the original component
contains the transitions \(a_i\xrightarrow{\,x_i\mid y\,}b\), \(i=1,2\).
Thus \(a_i\cdot x_i=b,
\,
x_i|_{a_i}=y.\)
The inverse-section identity gives \(x_1^{-1}|_b
=
y^{-1}
=
x_2^{-1}|_b.\) If \(x_1=x_2\), reversibility of the dual component would imply
\(a_1=a_2\), a contradiction. Hence \(x_1\neq x_2\). Since
\(\mathcal A\) is reduced, \(x_1^{-1}\) and \(x_2^{-1}\) represent
distinct elements of \(G(\mathcal A)\). Therefore
\(\operatorname{sec}_b\) is not injective and Theorem~\ref{theo:section-surjective-collision} now applies.
\end{proof}

In particular, if a finite reduced reversible invertible Mealy
automaton generates a non-abelian free group, then at least one
connected component of its dual is bireversible.

\subsection{A folding obstruction to non-reversibility}

The preceding criterion assumes reversibility. When the states form a
free basis, a transition collision can instead be detected by folding
the Schreier graph of a vertex stabilizer. Retain the notation
\(P_a\) and \(\sigma_a\) of Section~\ref{sec: section homomorphisms},
with \(W=X\).

\begin{lemma}
\label{lem:transition-collision-fold}
Let \(\mathcal A=(X,A,\delta,\lambda)\) be a finite invertible Mealy
automaton and assume that the natural map
\(\rho:F_X\to G(\mathcal A)\) is injective. If
\(x_1\neq x_2\) and \(x_1|_a=x_2|_a\) for some \(a\in A\), then
\(\ker(\sigma_a)\neq1\).
\end{lemma}

\begin{proof}
Identify \(G(\mathcal A)\) with \(F_X\) through \(\rho\), and consider the orbit 
\(O:=a\cdot G(\mathcal A)\), and put \(d:=|O|\), and \(r:=|X|\). Let \(\Gamma=\operatorname{Sch}(F_X,P_a,X)\) be the Schreier graph of
the action of \(F_X\) on \(a\cdot F_X\), based at \(a\). Then closed
paths at \(a\) are labelled precisely by the elements of \(P_a\), and
hence \(\pi_1(\Gamma,a)\simeq P_a.\)
Relabel each positive edge \(v\xrightarrow{\,x\,}v\cdot x\) of
\(\Gamma\) by \(x|_v\); its reverse edge is labelled by
\((x|_v)^{-1}\), since
\(x^{-1}|_{v\cdot x}=(x|_v)^{-1}\). By the product identity, the
label of a closed path based at \(a\) is the image under \(\sigma_a\)
of its original label. Hence folding this inverse labelled graph and
taking its core gives the Stallings graph of \(\sigma_a(P_a)\).

The two edges issued from \(a\) and originally labelled by \(x_1\)
and \(x_2\) have the same new label. If they have the same terminal
vertex, their fold decreases the rank. Otherwise the first fold
identifies their terminal vertices. The resulting graph has at most
\(d-1\) vertices, and further folds cannot increase this number. A
connected deterministic inverse \(X\)-labelled graph with \(M\)
vertices has at most \(rM\) geometric edges, and therefore rank at
most \(1+(r-1)M\). In either case,
\[
\rk(\sigma_a(P_a))<1+d(r-1)=\rk(P_a).
\]
Thus \(\sigma_a:P_a\to\sigma_a(P_a)\) is not injective.
\end{proof}
\begin{theorem}
\label{theo:free-basis-bireversible}
Let \(\mathcal A=(X,A,\delta,\lambda)\) be a finite invertible Mealy
automaton with \(|X|\geq2\). Assume that the natural epimorphism \(\rho:F_X\longrightarrow G(\mathcal A)\) is an isomorphism and that \(G(\mathcal A)\) acts transitively on \(A\).
Then \(\mathcal A\) is bireversible.
\end{theorem}
\begin{proof}
If \(\mathcal A\) were not reversible, there would exist \(a\in A\)
and \(x_1\neq x_2\) such that \(x_1|_a=x_2|_a\). By
Lemma~\ref{lem:transition-collision-fold}, \(\ker(\sigma_a)\neq1\).
Since the action on \(A\) is transitive,
Theorem~\ref{theo: local-global-criterion}, with \(n=1\), would imply
\(\ker(\rho)\neq1\), contradicting the injectivity of \(\rho\).  Since \(\rho\) is injective, \(\mathcal A\) is reduced, and
\(G(\mathcal A)\) is a non-abelian free group. Moreover,
\(\partial\mathcal A\) is connected because the action on \(A\) is
transitive. If its unique component were non-bireversible,
Corollary~\ref{cor:reversible-dual} would imply that
\(G(\mathcal A)\) is cyclic or not free, a contradiction. Hence
\(\mathcal A\) is bireversible.
\end{proof}

\subsection{The level-transitive binary case}
\label{subsec:binary-case}

We now remove both reversibility and the free-basis assumption over a
binary alphabet, under level-transitivity. We call a subgroup of
$\Aut(A^*)$ finite-state if all its elements are finite-state. The
first lemma determines the index of the image of a first-level section
homomorphism. Its proof is the only place where we use a growth-gap
theorem.

\begin{lemma}
\label{lem:section-image-index}
Let $G\leq\Aut(A^*)$ be a finitely generated finite-state,
level-transitive non-abelian free group. For $a\in A$, put
$K_a:=\Stab_G(a)$ and $L_a:=\theta_a(K_a)$. Then $\theta_a$ is
injective and $[G:L_a]=|A|$.
\end{lemma}

\begin{proof}
Let $W$ be a free basis of $G$. If $\ker(\theta_a)\neq1$, then the
corresponding lifted section homomorphism $\sigma_a$ has non-trivial
kernel. Since $G$ is transitive on $A$,
Theorem~\ref{theo: local-global-criterion}, with $n=1$, contradicts
that $W$ is a free basis. Thus $\theta_a$ is injective.
\\
Let $S=\{\,g|_u:g\in W\cup W^{-1},\ u\in A^*\,\}\setminus\{1\}$ be the set of non-trivial sections of the elements of
$W\cup W^{-1}$. Then $S$ is finite, symmetric, generates $G$, and is closed under non-trivial sections. Hence
$|\theta_a(k)|_S\leq |k|_S$ for every $k\in K_a$. Since $\theta_a$ is
injective,
\[
 |K_a\cap B_S(n)|\leq |L_a\cap B_S(n)|.
\]
Now $[G:K_a]=|A|$, so $K_a$ and $G$ have the same exponential growth
rate with respect to $S$. Hence $L_a$ has the same exponential growth
rate as $G$. Since \(L_a=\theta_a(K_a)\) is finitely generated, it is quasiconvex in
the free group \(G\). If \([G:L_a]=\infty\), the growth-gap theorem \cite[Theorem~1.1]{DahmaniFuterWise}gives a strictly smaller exponential growth rate
for \(L_a\) than for \(G\), contradicting the preceding inequalities. Hence \([G:L_a]<\infty\). Put $d:=[G:L_a]$ and $r:=\rk(G)$. Since $K_a\simeq L_a$, Schreier's
formula gives
\[
 1+|A|(r-1)=\rk(K_a)=\rk(L_a)=1+d(r-1).
\]
As $r\geq2$, we obtain $d=|A|$.
\end{proof}
For the rest of the subsection, let $A=\{0,1\}$ and
$\sigma=(0\,1)$. For $G\leq\Aut(A^*)$, put
\[
 \mathcal D(G):=\langle g|_a:g\in G,\ a\in A\rangle.
\]

\begin{lemma}
\label{lem:binary-descent-step}
Let $G\leq\Aut(A^*)$ be a finitely generated finite-state,
level-transitive non-abelian free group. Suppose that
$h=(1,s)\sigma\in G$. Then $H:=\mathcal D(G)$ is a finite-state,
self-similar, level-transitive non-abelian free group, with
$[G:H]=2$ and $s\in H$.
\end{lemma}

\begin{proof}
Put \(K:=\Stab_G(0)=\Stab_G(1)\), where the equality follows since $A=\{0,1\}$, and let $H_i:=\theta_i(K)$ for $i\in\{0,1\}$. For $k\in K$, the section
identities give $(hkh^{-1})|_0=k|_1$. Since $hKh^{-1}=K$, we obtain
$H_0=H_1$; denote their common value by $H$.
Lemma~\ref{lem:section-image-index} gives $[G:H]=2$. Moreover,
$h^2=(s,s)\in K$, so $s=\theta_0(h^2)\in H$.
\\
Every element of $G$ belongs to $K$ or to $hK$. For $k\in K$ we have
$k|_0,k|_1\in H$, while $(hk)|_0=k|_1$ and
$(hk)|_1=sk|_0$. Since $s\in H$, every first-level section of an
element of $G$ belongs to $H$. Conversely, every element of $H=H_0$ occurs as a section at $0$ of an element of $K$. Hence
$\mathcal D(G)=H$. If $l=g|_a\in H$, then $l|_b=g|_{ab}\in H$, so $H$ is self-similar.
Given $u,v\in A^n$, choose $g\in G$ with $0u\cdot g=0v$. Then
$g\in K$ and $u\cdot g|_0=v$, with $g|_0\in H$. Thus $H$ is
level-transitive. It is finite-state as a subgroup of $G$, and
Schreier's formula gives $\rk(H)=2\rk(G)-1\geq3$.
\end{proof}

\begin{lemma}
\label{lem:binary-terminal-reduction}
Let $G\leq\Aut(A^*)$ be a finitely generated finite-state,
level-transitive non-abelian free group. If
$\operatorname{sec}_0:G\to G$ is not injective, then there exist
subgroups $L<H\leq G$ and an element
$h=(1,s)\sigma\in H\setminus L$ such that
$\mathcal D(H)=L$, $[H:L]=2$, and $s\in L$. Moreover, both $H$ and
$L$ are finite-state, self-similar, level-transitive non-abelian free
groups.
\end{lemma}

\begin{proof}
Since $\operatorname{sec}_0:G\to G$ is not injective, choose $x\neq y$ in $G$ with $x|_0=y|_0$. The elements $x$ and $y$
cannot induce the same permutation of $A$: otherwise
$xy^{-1}\in\Stab_G(0)$ and $(xy^{-1})|_0=1$, contradicting
Lemma~\ref{lem:section-image-index}. After interchanging $x$ and $y$,
assume that $x$ fixes $A$ and $y$ exchanges its two letters. Then
$h:=x^{-1}y=(1,s)\sigma$ for some $s\in G$. Set $G_0:=G$ and, as long as $h\in G_i$, define
$G_{i+1}:=\mathcal D(G_i)$. By
Lemma~\ref{lem:binary-descent-step}, each inclusion
$G_{i+1}<G_i$ has index $2$, and every $G_i$ is finite-state,
self-similar, level-transitive and non-abelian free. It remains to prove that the descent terminates. Let \(W\) be a finite generating set of \(G\) and put
\[
\Omega:=\{w|_u:w\in W\cup W^{-1},\ u\in A^*\}.
\]
Since \(G\) is finite-state, \(\Omega\) is finite; moreover it is closed
under sections and inverses. Suppose that \(Q=\langle T\rangle\) with
\(T\subseteq\Omega\). If
\(q=t_1^{\varepsilon_1}\cdots t_\ell^{\varepsilon_\ell}\in Q\), then,
for every \(a\in A\), the product identity expresses \(q|_a\) as a
product of first-level sections of elements of \(T\cup T^{-1}\).
Therefore
\[
D(Q)=
\left\langle
t|_a:t\in T\cup T^{-1},\ a\in A
\right\rangle.
\]
Since \(\Omega\) is closed under sections and inverses, \(D(Q)\) is
generated by a subset of \(\Omega\). Since \(W\subseteq\Omega\),
inductively every \(G_i\) is generated by a subset of the fixed finite set \(\Omega\). Hence the set of possible subgroups \(G_i\) is finite. As long as \(h\in G_i\), Lemma~\ref{lem:binary-descent-step} gives
\([G_i:G_{i+1}]=2\), so \(G_{i+1}<G_i\). Thus the strictly decreasing chain cannot continue indefinitely, and there exists \(m\) such that \(h\in G_m\setminus G_{m+1}\). Since the chain is strictly decreasing as long as $h\in G_i$, it must terminate. Thus, for some $m$, $h\in G_m\setminus G_{m+1}$. Set
$H:=G_m$ and $L:=G_{m+1}$. Then
$L=\mathcal D(H)$, $[H:L]=2$, and $s\in L$ by
Lemma~\ref{lem:binary-descent-step}; the same lemma gives all the
remaining assertions.
\end{proof}

\begin{lemma}
\label{lem:binary-terminal-bijective}
Under the conclusions of Lemma~\ref{lem:binary-terminal-reduction}, the maps $\operatorname{sec}_0,\operatorname{sec}_1:L\to L$ are
bijective.
\end{lemma}

\begin{proof}
Put $K:=\Stab_H(0)$ and $M:=K\cap L$. Since
$[H:K]=[H:L]=2$, both $K$ and $L$ are normal in $H$. Moreover
$K\neq L$, since $K$ fixes the first level whereas $L$ is
level-transitive. Hence $KL=H$. Therefore
\[
 K/(K\cap L)\simeq H/L,
 \qquad
 L/(K\cap L)\simeq H/K,
\]
and consequently $[K:M]=[L:M]=2$. By Lemmas~\ref{lem:section-image-index} and
\ref{lem:binary-descent-step}, both maps
$\theta_0,\theta_1:K\to L$ are isomorphisms. Put
$N:=\theta_0(M)$. Since $L\trianglelefteq H$ and $hKh^{-1}=K$, conjugation by $h$
preserves $M=K\cap L$. For $k\in K$ we have
$\theta_0(hkh^{-1})=\theta_1(k)$. Hence, for $k\in M$,
$\theta_1(k)\in\theta_0(M)$. Since conjugation by $h$ is a bijection
of $M$, we obtain
\[
 \theta_1(M)=\theta_0(M)=:N.
\]
As $\theta_0:K\to L$ is an isomorphism and $[K:M]=2$, it follows that
$[L:N]=2$. Moreover, $h\notin L$ and $[H:L]=2$, so $h^2\in L$. Since
$h^2=(s,s)\in K$, we have $h^2\in M$, and therefore
$s=\theta_0(h^2)\in N$. Finally, $H=K\sqcup hK$. Intersecting with $L$ gives
\[
 L=M\sqcup(L\cap hK).
\]
Since \(L\triangleleft H\), \([H:L]=2\), and \(h\notin L\), the coset
\(h^{-1}L\) is \(H\setminus L\). Hence, for \(k\in K\),
\[
hk\in L
\iff k\in h^{-1}L
\iff k\notin L
\iff k\in K\setminus M..
\]
Hence
\[
 L=M\sqcup h(K\setminus M).
\]
On $M$, the maps $\operatorname{sec}_0$ and
$\operatorname{sec}_1$ are the restrictions of $\theta_0$ and
$\theta_1$, respectively. Since
$\theta_0(M)=\theta_1(M)=N$, both are bijections $M\to N$.
Now let $k\in K\setminus M$. Then
$(hk)|_0=k|_1$ and $(hk)|_1=sk|_0$. Since
$\theta_0,\theta_1:K\to L$ are bijections and both send $M$ onto
$N$, their restrictions send $K\setminus M$ bijectively onto
$L\setminus N$. Moreover $[L:N]=2$, so $N\trianglelefteq L$; since
$s\in N$, left multiplication by $s$ preserves $L\setminus N$.
Therefore both section maps send $h(K\setminus M)$ bijectively onto
$L\setminus N$.
Thus $\operatorname{sec}_0,\operatorname{sec}_1:L\to L$ are
bijective.
\end{proof}

\begin{lemma}
\label{lem:binary-terminal-obstruction}
Let $H\leq\Aut(A^*)$ be a free group and let $L\trianglelefteq H$
have index $2$. Assume that $L$ is finitely generated, finite-state,
self-similar, level-transitive, and section-bijective on $A^*$. Then
there is no element $h=(1,s)\sigma\in H\setminus L$ with $s\in L$.
\end{lemma}

\begin{proof}
Assume, towards a contradiction, that
$h=(1,s)\sigma\in H\setminus L$ with $s\in L$, as in the statement.
Since $L\trianglelefteq H$, conjugation by $h$ restricts to an
automorphism $\beta:L\to L$, given by $\beta(x)=hxh^{-1}$.

Put $d:=h^2$. Since $[H:L]=2$ and $h\notin L$, we have $d\in L$.
Moreover, $h=(1,s)\sigma$ gives $d=(s,s)$. Since $H$ is free and
$h\neq1$, it is torsion-free, and therefore $ d\neq 1$. Finally, $\beta(d)=d$ and, for every $x\in L$,
\[
 \beta^2(x)=dxd^{-1}.
\]

We first use section-bijectivity to obtain automorphisms of a finite
valency Cayley graph of $L$. Choose a finite symmetric generating set
$Q\subseteq L\setminus\{1\}$ closed under first-level sections. Such
a set exists because $L$ is finite-state: take all sections of a
finite symmetric generating set and remove $1$. Since the section
maps are injective, no non-trivial element has a trivial section, so
the resulting set is still closed under sections. Let $\Gamma$ be the underlying unlabelled graph of
$\operatorname{Cay}(L,Q)$. For $a\in A$, define
$\rho_a:L\to L$ by $\rho_a(g):=g|_a$. We claim that $\rho_a$ is a
graph automorphism of $\Gamma$ fixing $1$. Indeed, $\rho_a$ is bijective by section-bijectivity. If \(g\xrightarrow{q}gq\) is an edge of the Cayley graph, then
\[
 \rho_a(gq)=\rho_a(g)\,q|_{a\cdot g}.
\]
For fixed \(a\) and \(g\), the map
\[
 Q\longrightarrow Q,\qquad q\longmapsto q|_{a\cdot g},
\]
is injective by section-bijectivity and has image in \(Q\), since
\(Q\) is section-closed. As \(Q\) is finite, it is a permutation of
\(Q\). Hence \(\rho_a\) maps the neighbours of \(g\) bijectively onto
the neighbours of \(\rho_a(g)\). Since \(\rho_a:L\to L\) is bijective,
it is a graph automorphism of \(\Gamma\), fixing \(1\). Moreover $\rho_a(1)=1$. Consequently
$\gamma:=\rho_1^{-1}\rho_0\in\operatorname{Aut}(\Gamma)$ fixes $1$.
In particular,
\[
 |\gamma^n(g)|_Q=|g|_Q
 \qquad(g\in L,\ n\in\mathbb Z).
\]
Let $\varepsilon:L\to C_2$ be the homomorphism induced by the action
on the first level, so that $x$ acts on $A=\{0,1\}$ as
$\sigma^{\varepsilon(x)}$. We now compute $\gamma$. Let $x\in L$. Since
$h|_0=1$, $0\cdot h=1$, $h^{-1}|_1=1$, and
$h^{-1}|_0=s^{-1}$, the section identities give
\[
 \rho_0(\beta(x))
 =
 \begin{cases}
  \rho_1(x),&\varepsilon(x)=0,\\
  \rho_1(x)s^{-1},&\varepsilon(x)=1.
 \end{cases}
\]
Now $d^{-1}=(s^{-1},s^{-1})$, so this may be written as
$\rho_0(\beta(x))=\rho_1(xd^{-\varepsilon(x)})$. Hence
\[
 \gamma(\beta(x))=xd^{-\varepsilon(x)}.
\]
Conjugation by $h$ preserves the permutation induced on the binary
first level, so $\varepsilon(\beta(x))=\varepsilon(x)$. Replacing
$x$ by $\beta^{-1}(y)$ therefore gives
\begin{equation}
\label{eq:binary-gamma}
 \gamma(y)=\beta^{-1}(y)d^{-\varepsilon(y)}
 \qquad(y\in L).
\end{equation}
Since $\beta(d)=d$ and $\beta^2(x)=dxd^{-1}$, we have
$\beta^{-2}(x)=d^{-1}xd$. Applying
\eqref{eq:binary-gamma} twice and using \(\varepsilon(d)=0\) gives
\[
 \gamma^2(y)
 =
 d^{-1}y\,d^{1-2\varepsilon(y)}.
\]
Since $\varepsilon(\gamma(y))=\varepsilon(y)$, induction yields
\begin{equation}
\label{eq:binary-gamma-powers}
 \gamma^{2n}(y)
 =
 d^{-n}y\,d^{n(1-2\varepsilon(y))}
 \qquad(n\geq0).
\end{equation}

Since $L$ is transitive on $A$, choose $y\in L$ with
$\varepsilon(y)=1$. Then
$\gamma^{2n}(y)=d^{-n}yd^{-n}$ for every $n\geq0$. The graph
automorphism $\gamma$ fixes $1$, hence preserves the distance from
$1$. Thus all these elements lie in the finite sphere of radius
$|y|_Q$ in $\Gamma$. There exist therefore $m>n$ such that
$d^{-m}yd^{-m}=d^{-n}yd^{-n}$. Setting $k:=m-n>0$ and simplifying,
we obtain
\[
 y^{-1}d^{-k}y=d^{k}.
\]
Let $r$ be the primitive root of $d$, so that $\langle r\rangle$ is
the unique maximal cyclic subgroup containing $d^k$. We have $y^{-1}\langle d^k\rangle y=\langle d^k\rangle$. Hence
$y^{-1}\langle r\rangle y$ is a maximal cyclic subgroup containing
$d^k$. By uniqueness,
$y^{-1}\langle r\rangle y=\langle r\rangle$. Thus $y$ normalizes
$\langle r\rangle$. Since maximal cyclic subgroups of free groups are self-normalizing, $y\in\langle r\rangle$. Since $y\in\langle r\rangle$ and $d\in\langle r\rangle$, we get $d^{2k}=1$. Since $H$ is
torsion-free, $d=1$, a contradiction.
\end{proof}

\begin{theorem}
\label{theo:binary-section-bijective}
Let $G\leq\Aut(\{0,1\}^*)$ be a finitely generated finite-state,
level-transitive non-abelian free group. Then
$\operatorname{sec}_u:G\to G$ is bijective for every
$u\in\{0,1\}^*$.
\end{theorem}

\begin{proof}
If $\operatorname{sec}_0$ were not injective,
Lemma~\ref{lem:binary-terminal-reduction} would give
$L<H\leq G$ and $h=(1,s)\sigma\in H\setminus L$.
Lemma~\ref{lem:binary-terminal-bijective} makes $L$
section-bijective, contradicting
Lemma~\ref{lem:binary-terminal-obstruction}. Thus
$\operatorname{sec}_0$ is injective. Interchanging $0$ and $1$ gives
the injectivity of $\operatorname{sec}_1$.

For surjectivity, put $K:=\Stab_G(0)$ and
$L_i:=\theta_i(K)$ for $i\in\{0,1\}$. By
Lemma~\ref{lem:section-image-index}, both $L_0$ and $L_1$ have index
$2$ in $G$. Choose $t\notin K$. Since $K$ and $tK$ are the two cosets
of $K$, we have
\[
 \operatorname{sec}_0(K)=L_0,
 \qquad
 \operatorname{sec}_0(tK)=t|_0L_1.
\]
Injectivity of $\operatorname{sec}_0$ gives
$L_0\cap t|_0L_1=\varnothing$. If \(L_0\neq L_1\), then \(L_0L_1=G\), since distinct subgroups of
index \(2\) are normal and generate \(G\). Thus, for every \(g\in G\),
we may write \(g=l_0l_1\) with \(l_i\in L_i\), so
\(l_0=gl_1^{-1}\in gL_1\cap L_0\). Hence every coset of \(L_1\)
meets \(L_0\), contradicting
\(L_0\cap t|_0L_1=\varnothing\). Thus $L_0=L_1$, and $t|_0L_1$ is the other coset of $L_0$.
Therefore $\operatorname{sec}_0(G)=G$. The same argument applies to $\operatorname{sec}_1$. Finally, if $u=a_1\cdots a_n$, then
$\operatorname{sec}_u=
\operatorname{sec}_{a_n}\circ\cdots\circ\operatorname{sec}_{a_1}$,
so $\operatorname{sec}_u$ is bijective.
\end{proof}
\begin{cor}
\label{cor:binary-bireversible}
Let $\mathcal A$ be a finite reduced invertible Mealy automaton over a
binary alphabet.  If $G(\mathcal A)$ is a non-abelian free group and
acts level-transitively, then $\mathcal A$ is bireversible.
\end{cor}

\begin{proof}
Set $G:=G(\mathcal A)$.  By
Theorem~\ref{theo:binary-section-bijective}, the maps
$\operatorname{sec}_0$ and $\operatorname{sec}_1$ are injective on
$G$.  If $x_1|_a=x_2|_a$ for two states of $\mathcal A$, then
$x_1=x_2$ in $G$, and reducedness gives equality of the states. Thus $\mathcal A$ is reversible. Suppose that the inverse automaton is not reversible.  Then the
original automaton has transitions \( x_i\xrightarrow{\,a_i\mid b\,}y, i=1,2\) with $x_1\neq x_2$.  The inverse-section identity gives
\[
 x_1^{-1}|_b=y^{-1}=x_2^{-1}|_b.
\]
Injectivity of $\operatorname{sec}_b$ gives $x_1=x_2$ in $G$, and
reducedness gives equality of the states, a contradiction.  Therefore
the inverse automaton is reversible as well.
\end{proof}
\begin{prob}\label{problema}
Let $G\leq\Aut(A^*)$ be a self-similar non-abelian free group. Must
there exist a $G$-orbit $B\subseteq A$ such that
$\operatorname{sec}_u:G\to G$ is bijective for every $u\in B^*$?
Corollary~\ref{cor:intrinsic-bijective-orbit} gives an affirmative
answer when the action is section-surjective on every first-level
orbit, while Theorem~\ref{theo:binary-section-bijective} gives one for
finite-state level-transitive actions on a binary tree.
\end{prob}
The corresponding automata-theoretic question is the following.
\begin{prob}
If a finite reduced invertible Mealy automaton generates a non-abelian
free group acting transitively on the first level, must its dual
contain a bireversible connected component?
Corollary~\ref{cor:reversible-dual} gives an affirmative answer when
the automaton is reversible, while
Theorem~\ref{theo:free-basis-bireversible} gives one when the states
form a free basis. Corollary~\ref{cor:binary-bireversible} settles the
binary case under level-transitivity. The first-level-transitive case
with a redundant state generating set remains open in general; under
level-transitivity, the corresponding redundant non-reversible case
on alphabets of size at least three also remains open.
\end{prob}

\begin{rem}
\label{rem:binary-square-complex}
Under the hypotheses of Corollary~\ref{cor:binary-bireversible}, let
$\Delta_{\mathcal A}$ denote the square complex associated with
$\mathcal A$ as in \cite{BoKi}. The corollary gives that
$\mathcal A$ is bireversible, while $G(\mathcal A)$ is infinite
because it is a non-abelian free group. Since the alphabet is binary,
\cite[Theorem~11]{BoKi} implies that
$\pi_1(\Delta_{\mathcal A})$ is not residually finite. On the other
hand, $G(\mathcal A)$ is residually finite, being a free group.
Thus $G(\mathcal A)$ is residually finite whereas
$\pi_1(\Delta_{\mathcal A})$ is not.
\end{rem}

The two obstructions developed in this paper describe complementary regimes for
the dynamics of sections. In the absorbing regime, where iterated sections
collapse to the identity, Theorem~A rules out freeness, while
Proposition~\ref{prop:activity-non-expansion} shows that the level Schreier graphs have
vanishing Cheeger constants, so that Kazhdan's property~(T) fails as well
(Corollary~\ref{cor:coaccessible-no-T}). Theorems~B, C and D concern the
opposite phenomenon: under the corresponding hypotheses, non-abelian freeness
forces bireversibility, or at least the presence of a bireversible component.
Bireversible automata are precisely where both types of algebraic behaviour are
known to occur: the square-complex method of \cite{Mozes} produces bireversible
automata generating non-abelian free groups, as well as bireversible automata
whose groups have property~(T). Thus both known free and property-\((T)\)
examples lie outside the absorbing regime. For free groups, the results above
further suggest that bijective section dynamics should occur on some
first-level orbit; Problem~\ref{problema} is the precise formulation of this
expectation.

%%%%%%%%%%%%%%%%%%%%%%%%%%%%%%%%%%%%%%%%%%%%%%%%%%%

\section*{Acknowledgements}
The authors thank one referee for encouraging a substantial reframing
of an early version of the paper in more standard terms, and another
for encouraging them to address an open problem stated in that
version, a suggestion that eventually led to several of the additional
results presented here. The authors are members of the National Research Group
GNSAGA (Gruppo Nazionale per le Strutture Algebriche, Geometriche e le
loro Applicazioni) of INdAM.

\end{document}